\documentclass[3p]{elsarticle}

\usepackage{amsmath,amsfonts,amssymb,amsbsy,amscd}
\usepackage{graphicx}
\usepackage{mathrsfs,array}
\usepackage{verbatim} 
\usepackage{subfig}

\usepackage[section,subsection,subsubsection]{extraplaceins}

\def\B#1{\mbox{\boldmath{$#1$}}}

\newcommand{\bu}{\B{u}}

\newcommand{\bx}{\B{x}}

\newcommand{\onehalf}{\mbox{$\frac{1}{2}$}}

\def\be{\begin{equation}}
\def\ee{\end{equation}}
\def\ba{\begin{array}}
\def\ea{\end{array}}
\def\bea{\begin{eqnarray}}
\def\eea{\end{eqnarray}}
\def\beas{\begin{eqnarray*}}
\def\eeas{\end{eqnarray*}}
\newcommand{\bseq}{\begin{subequations}}
\newcommand{\eseq}{\end{subequations}}

\begin{document}

\begin{frontmatter}

\title{Monotone level-sets on arbitrary meshes without redistancing}

\author{Ido Akkerman}
\address{ Mechanical, Maritime and Materials Engineering Department  \\ 
 Delft University of Technology}
\date{\today}

\begin{abstract}
In this paper we present approaches that address two issues
that can occur when the level-set method is used to simulate two-fluid 
flows in engineering practice.

The first issue concerns regularizing the Heaviside function on arbitrary meshes.
We show that the regularized Heaviside function can be non-smooth on non-uniform meshes. 
 Alternative regularizing definitions that are indeed smooth and monotonic, are introduced.  
 These new definitions lead to smooth Heaviside functions by taking the changing local meshsize into account.

The second issue is the computational cost and fragility caused by the necessity of redistancing the level-set field.
In \cite{AkBaBeFaKe12,Akkerman2012} it is shown that strongly coupling the level-set 
convection with the flow solver provides robustness and potentially efficiency and accuracy advantages.
The next step would be to include redistancing within the strong coupling part of the algorithm.
The computational cost  of current redistancing procedure prohibit this.
Four alternative approaches for circumventing the expensive redistancing  step are proposed.
This  should facilitate a fully coupled level-set approach.

Some benchmark cases  demonstrate the efficacy of the proposed approaches.
These includes the standard test case of the vortex in a box. Based on these results 
the most favourable redistancing approach is selected. 
\end{abstract}


\end{frontmatter}

\section{Introduction}
\label{sec:intro}

Level-sets are a very powerful approach for solving interface problems \cite{Osher_Sethian_88}.   
By disconnecting the interface description from the underlying mesh, topology changes of the interface are handled with ease.
For numerical reasons, such as numerical quadrature or finite differencing, the interface is often 
given a finite mesh-dependent width. This is done by adopting a so-called regularized Heaviside function, where the transition 
from zero to unity does not occur instantaneously but over a finite band. 
This regularization, however, causes two problems. On irregular meshes -- which might occur in the
 analysis of engineering artefacts  -- the smoothed Heaviside can become non-smooth or even non-monotone. 
 Additionally, the level-set is required to be a distance function in order to control the thickness of the interface. 
Due to its evolution the level-set  evolves requires  to be redistanced, that is 
the distance property needs to be actively repaired.

In \cite{AkBaBeFaKe12,Akkerman2012} the level-set method is used to model a rigid-body floating on a 
water surface. This water surface is handled using a level-set. In this paper we adopted a strongly coupled modeling approach. 
Meaning that both the flow; interface evolution; body motion and mesh deformation are solved 
simultaneously.
In agreement with \cite{Felippa01,TT_ARCH,Walhorn2005} -- and other FSI literature-- this approach results in a robust solution strategy. 
In a pure two-fluid problem, i.e. without a floating object, the strongly coupled approach proved to be more robust than traditional approaches
where interface evolution is solved after the flow problem is solved.

However, the redistancing of the level-set was not included in the strongly coupled part of the solver, this was done at the end of each time step.
Due to the redistribution of mass  associated with redistancing, momentum and energy conservation are difficult to control. 
Certainly, the energy errors are troublesome as these might trigger instabilities in the solver. 
Therefore, there is a desire to include the redistancing step in the strongly coupled part of the method. 
This should allow methods that are either energy conservative 
or guaranteed to be energy dissipative. This is envisioned to lead to more robust  and  more accurate methods.

In this paper we present four alternative approaches which lead to an approximate distance field without directly solving a redistancing problem.
These approaches are formulated in such a way that they potentially solve the non-smoothness on irregular grids.

The paper is structured as follows,
in section \ref{sec:LS} provides a short introduction to the level-set method and the regularized Heaviside function.
In section \ref{sec:Mono} we show that a naive definition of the regularized  Heaviside can actually be non-smooth.
A potential solution to this problem is presented.

In section \ref{sec:noRD} we build on the result from the previous section and introduce four alternative approaches 
to redistance the level-set field.  All these approaches circumvent  solving the difficult non-linear Eikonal equation. The Eikonal equation is translated to simple projection problems that achieve the same:  a distance field suitable to define a smooth Heaviside.

In section \ref{sec:NumAlg} we give a short description of the numerical formulations and the finite element and isogeometric discretizations employed to solve the convection and redistancing problems in the next section. In this section, 
section \ref{sec:NumExp},  we first solve a simple redistancing problem on irregular grids and then 
solve the vortex in a box problem using different discretizations and different redistancing approaches.
Using these test cases we identify the best method, and the appropriate numerical parameters.
For this method a mesh convergence study  is performed and suitability in three dimensions is tested.

In section \ref{sec:conc} we conclude and sketch a perspective for further use of the proposed methods

\section{The level set method}
\label{sec:LS}

The level-set method has since its origin \cite{Osher_Sethian_88} been applied to numerous problems involving surfaces, interfaces
and shapes. 
Its  application can be found in a wide range of areas  from image processing, computer graphics, 
topology and shape optimisation and simulation of physics problems at interfaces such as crystal growth, or two-fluid flow.  

See, the review papers \cite{SethianAnnRev, Osher2001463} or books \cite{Sethian99,OsherBook} for a broad picture of the applications and methods available.
In this paper we  focus on level-set  in combination with variational methods such as finite element and isogeometric analysis \cite{HuCoBa04} with the final goal of applying it to two-fluid flow. The two-fluid problem  is a physical problem where the 
interface between the fluids is unknown and physical parameters are discontinuous across the interface.
Adoption of a level-set allows for easy handling of interface topology changes and regularisation  of the discontinuity.

\subsection{Level-sets for two-fluid problems}

In the level-set method a surface of lower dimension, 
such as an interface between two distinct materials is indirectly 
parameterized employing a globally defined function. This function is denoted as $\phi$, and defines a surface as follows,
\begin{align}
\Gamma_i = \{ \bx \in \Omega : \phi = 0 \}.
\end{align}
This automatically leads to the following distinct subdomains,
\begin{align}
\Omega^- =& \{ \bx \in \Omega : \phi < 0 \} ,\nonumber \\
\Omega^+ =& \{ \bx \in \Omega : \phi > 0 \},
\end{align}
which allows
the prescription of different physical parameters in each subdomain. For instance, the density
\begin{align}
\rho = 
 \left \{
 \begin{array}{ccc}
\rho_0 &\text{if}& \bx \in \Omega^- , \\
\rho_1 &\text{if}& \bx \in \Omega^+ .
 \end{array}
\right .
\end{align}

Alternatively, the Heaviside function
\begin{align}\label{eq:sharpH} 
 H(\phi) =
 \left \{
 \begin{array}{ccc}
 0            &\text{if}& \phi<0, \\
 \onehalf &\text{if}& \phi = 0, \\
 1            &\text{if}& \phi >0,
 \end{array}
\right .
\end{align}
can be used for a  convex interpolation to define a density
\begin{align}
\rho = \rho_0 (1-H(\phi)) + \rho_1 H(\phi).
\end{align}
This has the advantage of automatically handling the interface itself in a natural way.

\subsection{Regularized Heaviside}

In numerical methods the sharp interfaces defined in the previous section can lead to problems. 
For instance this appears when determining a mass matrix, $M_{ab}$ is approximated by quadrature as,
\begin{align}
M_{ab} = \int_{\Omega} \rho N_a N_b d\Omega  \approx \sum_{i=1..n_{ip}} \rho(\bx_i) N_a(\bx_i) N_b(\bx_i) w_i.
\end{align}
A sudden change of the density leads to a very bad approximation of the intended integral.  

To alleviate this problem  the sharp Heaviside function, 
defined in eq (\ref{eq:sharpH}), is replaced by a regularized Heaviside function. This regularized Heaviside function is often defined as
\begin{align}\label{eq:smoothH1} 
 H_{\epsilon}(\phi) =
 \left \{
 \begin{array}{ccc}
 0            &\text{if}& \phi< - \epsilon, \\
 \onehalf(1+\sin (\frac{\pi \phi}{2 \epsilon})) &\text{if}& \phi = 0, \\
 1            &\text{if}& \phi >\epsilon,
 \end{array}
\right .
\end{align}
where $\epsilon$ is the smoothing distance. Instead of an instantaneous switch from 0 to 1, this switch is 
spread over a finite layer around the interface. To have strict control over 
the width of this interface layer, we require $\phi$, to be a 
signed distance function. This means it needs to satisfy
\begin{align}
\| \nabla\phi \| = 1, 
\end{align}
which is the Eikonal equation.
As the regularizing is introduced to deal with numerical issues, such as quadrature, 
it is natural to specify the finite interface layer in terms of mesh size $h$ as
\begin{align}
\epsilon = \alpha h. 
\end{align}
Here $\alpha$ is an $O(1)$ parameter. The 
meshsize can be defined unambiguously for structured equidistant meshes. 
However, on arbitrary meshes this is not always straightforward.
In  \cite{ABKF11,KAFB10,AkBaBeFaKe12} we employed the meshsize
\begin{align}\label{eq:h1} 
 h  =& \frac{\|\nabla \phi \|} {\sqrt{\nabla \phi \cdot \B{G} \nabla \phi}},
 \end{align}
 where $\B{G}$ is the metric-tensor
\begin{align} \label{eq:metric} 
  \B{G}& =\left (  \frac{\partial \xi}{\partial  \bx} \right)^T \frac{\partial \xi}{\partial  \bx}
\end{align}
where $\bx$ is the physical space coordinate and $\xi$ is the coordinate in parametric space pertaining to the reference element.

This definition of $h$ incorporates the desired directional information. Effectively, 
a length scale is extracted from the metric-tensor in the direction $\frac{\nabla \phi}{\|\nabla \phi \|}$ .

\section{Monotonicity on arbitrary meshes}
\label{sec:Mono}

In this section we further discuss the smoothing of the Heaviside  function on  arbitrary meshes.
For this exposition it is useful to slightly rewrite the regularized Heaviside function,
\begin{align}\label{eq:smoothH2} 
 \hat{H}(\hat{\phi}) =
 \left \{
 \begin{array}{ccc}
 0            &\text{if}& \hat{\phi}< -\alpha,\\
  \onehalf (1+\sin(\frac{\pi}{2} \frac{\hat{\phi}}{\alpha} )) &\text{if}& \hat{\phi}= 0, \\
 1            &\text{if}& \hat{\phi} >\alpha,
 \end{array}
\right .
\end{align}
where 
\begin{align}\label{eq:naiveScalling} 
\hat{\phi} = \frac{\phi}{h}
\end{align}
is the scaled distance. In  other words,  if $\phi$ is the {\it actual} distance to the
 interface,  expressed in for instance $mm$ or $m$,  then $\hat{\phi}$ can be thought of as that same distance but 
 expressed in {\it number of elements}, that is a multiple of a typical element length. However, since this scaling occurs locally, the variation
  of $h$ along the path - from interface to the point under consideration  - is not taken into account. 
  Therefore, this rescaled 
  distance is only an effective estimate if $h$  varies only mildly (or not at all) across the interface.

\begin{figure}[!ht]
\centering
\subfloat[][Scaled distance on a uniform mesh]
{
\input{smooth/uni_phi.tex}
\label{uni_phi}
}
\subfloat[][Smooth Heaviside on a uniform mesh]
{
\input{smooth/uni_H.tex}
\label{uni_H}
}
\\
\subfloat[][Scaled distance on a non-uniform mesh]
{
\input{smooth/arb1_phi.tex}
\label{arb1_phi}
}
\subfloat[][Smooth Heaviside on a non-uniform mesh]
{
\input{smooth/arb1_H.tex}
\label{arb1_H}
}
\caption{ Original scaled distance and regularized Heaviside on a uniform and non-uniform mesh.
}
\label{fig:org_LS}
\end{figure}

Figure \ref{fig:org_LS} visualizes the scaled distance and the resulting regularized Heaviside function
on an uniform and non-uniform mesh, respectively.  The element distribution is indicated in the 
two left figures.   
The top plots show the results on a uniform mesh. 
In this case the scaled distance is a straight line and the regularized Heaviside function is indeed smooth.
The bottom plots show the results obtained on a non-uniform mesh. Here 
the scaled distance  and more importantly the regularized Heaviside function show a saw-tooth behaviour.
This defeats the purpose of the regularization. 
Although, one could argue that the quadrature on each element behaves well as the jumps occur at element interfaces. 
This is indeed the case, however, this will result in an inconsistent numerical method unless the density jumps at the interface are handled correctly, i.e. analogous to in the  discontinous  Galerkin method \cite{DGLS}.

\subsection{Path dependent rescaling}

In order to solve the problems indicated in the previous section. 
We should translate eq (\ref{eq:naiveScalling}) to a pure local relation as
\begin{align}\label{eq:pathScalling0} 
d\hat{\phi} = \frac{d\phi}{h}.
\end{align}
The total scaled distance can be found by path integration,
\begin{align}\label{eq:pathScalling1} 
\hat{\phi} = \int \frac{d\phi}{h}
\end{align}
which is not trivial to compute in practice.

Instead we take another route to obtain $\hat{\phi}$.  As an alternative to eq (\ref{eq:pathScalling0}) we can also state,
\begin{align}
\nabla \hat{\phi} = \frac{\nabla \phi}{h}.
\end{align}
This immediately results in an alternative Eikonal equation for the scaled distance:
\begin{align}\label{eq:scaledEikonal} 
\| \nabla \hat{\phi} \| = \frac{\| \nabla \phi \| }{h} = \frac{1}{h},  
\end{align}
where we used $\| \nabla \phi \|=1$. In principle any reasonable selection of $h$ will work.

However, inspired on the definition for $h$ given in eq (\ref{eq:h1}), we can define $h$ in such a way that    
eq (\ref{eq:scaledEikonal}) gets replaced with an attractive alternative Eikonal equation.
To achieve this we choose 
\begin{align}\label{eq:h2} 
 h  =& \frac {\sqrt{\nabla\hat{\phi} \cdot \B{G}^{-1} \nabla \hat{\phi}}}{\|\nabla \hat{\phi} \|},
 \end{align}
 where
 \begin{align}\label{eq:inv-metric} 
  \B{G}^{-1}& =  \frac{\partial \bx}{\partial \xi} \left ( \frac{\partial \bx}{\partial  \xi}\right)^T
\end{align}
is the inverse  of the metric-tensor defined in eq (\ref{eq:metric}). 
The same rational behind eq (\ref{eq:h1}) is employed eq (\ref{eq:h2}).  
We use the direction $\frac{\nabla \hat{\phi}}{\|\nabla \hat{\phi}\| }$ 
to extract the relevant mesh size from  the inverse metric tensor. 
The metric tensor depends on the rotation, scaling and distortion of the element and are naturally included in the mesh size definition. 

Defining $\nabla_{\xi}$ as the gradient in the reference direction, we can use the chain rule to arrive at
\begin{align}
\nabla_{\xi} \hat{\phi} \equiv
 \frac{\partial \hat{\phi}}{\partial \xi} =
 \frac{\partial \hat{\phi}}{\partial \bx} 
 \frac{\partial \bx}{\partial \xi}  = \nabla \hat{\phi}   \frac{\partial \bx}{\partial \xi},
\end{align}
This is substituted into eq (\ref{eq:h2}) to obtain 
\begin{align}
 h  =&\frac {\|\nabla_{\xi} \hat{\phi} \|} {\|\nabla\hat{\phi} \|}.
\end{align}
Combining this with the scaled Eikonal equation  (\ref{eq:scaledEikonal}) we get
\begin{align}\label{eq:meshEikonal} 
\|\nabla_{\xi} \hat{\phi} \|=1,
\end{align}
which is the Eikonal equation in terms of the parameteric space. 
This equation states  that $\hat{\phi}$ is a distance in terms of mesh lengths, 
which was exactly what we were looking for in the first place.

\begin{figure}[!ht]
\centering
\subfloat[][Corrected scaled distance.]
{
\input{smooth/arb1_phi2.tex}
\label{arb1_phi2}
}
\subfloat[][Corrected regularized  Heaviside.]
{
\input{smooth/arb1_H2.tex}
\label{arb1_H2}
}
\caption{Modified scaled distance and regularized  Heaviside on  non-uniform mesh.
\label{fig:New_LS}
}

\end{figure}

Figure \ref{fig:New_LS} shows the resulting scaled distance and regularized Heaviside when the alternative scaling is used.
The mesh is the same as in  figure \ref{fig:org_LS}. Both the distance and Heaviside function do not show a saw-tooth behaviour and the the Heaviside function is actually smooth, as intended. 
The scaled distance is indeed not a straight line, which is caused by the metric defined by the non-uniform mesh.

\section{Circumventing explicit redistancing}
\label{sec:noRD}

In the previous section we have addressed the issue of monotone smoothing of the Heaviside function. 
An equally important issue is that of the distance property of $\phi$,  or better of $\hat{\phi}$. 
The level-set time-evolution is typically given by
\begin{align}
\label{eq:convection}
\frac{\partial \phi}{\partial t} + \bu\cdot \nabla\phi = 0 
\end{align}
or something equivalent.  It is not guaranteed whether its distance property is maintained during the simulation.  
In order to keep explicit control of the interface layer thickness, we need to restore the distance property. 
Even if $\phi$ would remain a distance function, the path dependent scaling of the distance function
still requires us to solve an Eikonal equation, namely eq (\ref{eq:meshEikonal}), 
\begin{align}
\|\nabla_{\xi} \hat{\phi} \|=1 \qquad \text{on} ~~ \Omega^- \cup \Omega^+.
\end{align}
This equation is augmented with the following embedded boundary condition,
\begin{align}
\hat{\phi} =0  \qquad \text{on} ~~ \Gamma_i= \{ \bx \in \Omega : \phi = 0 \},
\end{align}
which forces the interface to remain fixed in the scaling/redistancing step.

Many different approaches to solve the Eikonal equation have been developed, 
see for instance \cite{Sethian_01} and references therein.
Most methods introduce a  pseudotime $\tau$ and march in time using
\begin{align}
\frac{\partial  \hat{\phi}}{\partial \tau} +  \B{a} \cdot \nabla \hat{\phi}= S,
\end{align}
where $\B{a} $ is an effective convection velocity and $S$ is a source. Both $\B{a} $ and  $S$ are carefully 
chosen to guarantee convergence towards the desired solution. 

Remark: to guarantee a steady state solution  the source term needs to include a term penalizing
the deviation of the interface, as shown in \cite{AkBaBeFaKe12}.
This issue is often circumvented by limiting the number of pseudotime steps.

This redistancing procedure has several drawbacks.
It is quite costly, even when limited in space (only redistance near the interface) 
or time (only redistance every so many  timesteps). Additionally, careful attention needs to
 be paid to the interplay of different  parameters associated with the time-marching 
 algorithm.  This makes the procedure somewhat fragile. Last but not least, it is difficult to 
 find the right balance between redistancing on the one hand and maintaining the interface on the other hand. 
Therefore, it is preferred   to  circumvent explicit redistancing.

\subsection{General idea}
In the following sections we present four alternative procedures to achieve this.  
All four alternatives hinge on the idea that 
 we  relate the unknown $\hat{\phi}$ with the known $\phi$ by
\begin{align}\label{eq:hat-phi}
\hat{\phi} = \phi \epsilon,
\end{align}
where $\epsilon$ is an unknown multiplicative scaling, which is required to be strictly positive. 
Clearly,
$\hat{\phi}$  vanishes when $\phi$ vanishes. This means that the interface location is trivially conserved.
There is no need to balance redistancing with maintaining the interface. 

The multiplicative scaling can be found by enforcing the Eikonal equation in parametric space (\ref{eq:meshEikonal}), 
\begin{align}
\|\nabla_{\xi} \hat{\phi}\|=1 \qquad \text{on}  \qquad \Omega^- \cup \Omega^+
\end{align}
 and subsequently  applying the chain rule.
As we are only interested in what happens near the interface we assume $\phi\approx 0 $
which results in the following relation for the multiplicative scaling,
\begin{align}
\epsilon = \frac{1}{\|\nabla_{\xi} \phi\|}
\end{align}
In the following section four different  alternative approaches based on this general idea are presented.

\subsection{First alternative: Direct  redistancing}

The first option is straightforward pointwise scaling:
\begin{align}
\hat{\phi} = \frac{\phi} {\|\nabla_{\xi} \phi \|}
\end{align}
This has the risk of resulting in divisions by zero.
Additionally, on standard $C_0$ continuous discretizations this will likely lead to
a discontinuous $\hat{\phi}$  field  and  consequently a discontinuous $H(\hat{\phi} )$ field.
Therefore this does not seem to be an attractive option.

\subsection{Second alternative: Projected redistancing}

We can obtain a continuous distance field by projecting the discontinuous $\hat{\phi}$  
of the previous section on a $C_0$ discretization.
 This can be done by solving the following projection problem
\begin{align}
\hat{\phi} - \nabla_{\xi}  \kappa_d \nabla_{\xi} \hat{\phi} =\frac{\phi} {\|\nabla_{\xi} \phi \|},
\end{align}
or by solving the following equivalent weak form:
\begin{center}
\fbox{\begin{minipage}{10cm}
{\it  Find $\hat{\phi} \in \mathcal{V}^h$ such that:}
\begin{align}
(w, \hat{\phi})_{\Omega} 
+ (\nabla_{\xi} w, \kappa_d \nabla_{\xi} \hat{\phi}) 
= \left (w,\frac{\phi} {\|\nabla_{\xi} \phi \|} \right )_{\Omega}  \qquad \forall w \in \mathcal{V}^h.
\end{align}
\end{minipage}}
\end{center}
Here, and in the remainder of this paper,  $\kappa_d$ is a smoothing parameter and $\mathcal{V}^h$ is an appropriate discrete space.
Due to the projection and the associated discretization errors, the interface location remains 
at the same positionin an approximate sense.
This means the interface is not necessarily  maintained during the redistancing step.

\subsection{Third alternative: Projected scaling}

By projecting the scaling coefficient on a continuous discretization the properties of continuity and interface conservation can be achieved simultaneously. This can achieved by either the following projection,
 \begin{align}
\epsilon - \nabla_{\xi}  \kappa_d \nabla_{\xi} \epsilon,
= \frac{1} {\|\nabla_{\xi} \phi \|} 
\end{align}
or by solving the equivalent weak form:
\begin{center}
\fbox{\begin{minipage}{10cm}
{\it  Find $\epsilon \in \mathcal{V}^h$ such that:}
 \begin{align}
(w, \epsilon)_{\Omega} 
+ (\nabla_{\xi} w, \kappa_d \nabla_{\xi} \epsilon ) 
= \left (w,\frac{1} {\|\nabla_{\xi} \phi \|} \right )_{\Omega}  \qquad \forall w \in \mathcal{V}^h.
\end{align}
\end{minipage}}
\end{center}
The sought after distance field can easily be computed as,
\begin{align}
\hat{\phi} = \phi \epsilon.
\end{align}

\subsection{Fourth alternative: Projected inverse scaling}

The last alternative approach is a small modification of the previous alternative.
Instead of projecting the multiplicative correction, we project its reciprocal value. 
This is achieved with the projection,
 \begin{align}
 \epsilon - \nabla_{\xi} \kappa_d \nabla_{\xi} \epsilon 
= \|\nabla_{\xi} \phi \| 
\end{align}
or by solving its corresponding weak form,
 \begin{center}
\fbox{\begin{minipage}{10cm}
{\it  Find $\epsilon \in \mathcal{V}^h$ such that:}
 \begin{align}
(w, \epsilon)_{\Omega} 
+ (\nabla_{\xi} w, \kappa_d \nabla_{\xi} \epsilon ) 
= \left (w,\|\nabla_{\xi} \phi \| \right )_{\Omega}  \qquad \forall w \in \mathcal{V}^h
\end{align}
\end{minipage}}
\end{center}
The sought after distance field is computed as
\begin{align}
\hat{\phi} = \frac{\phi }{\epsilon },
\end{align}
in this case.

\section{Numerical algorithm}
\label{sec:NumAlg}

In order to solve the convection problem given in eq (\ref{eq:convection}) we choose 
to use  finite elements and isogeometric analysis.
For completeness we describe the adopted weak formulation and associated numerical parameters. 
Next, we briefly discuss the isogeometric analysis shape functions that are used. 
Adoption of isogeometric analysis allows the selection of $C_1$ continuous discretizations.
This will turn out to be beneficial.

\subsection{SUPG formulation}

The convective problem from eq (\ref{eq:convection}) is solved using a SUPG formulation \cite{BroHug82}.
This problem is stated as follows:
\begin{center}
\fbox{\begin{minipage}{12cm}
{\it  Find $\phi \in \mathcal{V}^h$ such that:}
\begin{align}
(w+\tau \bu\cdot\nabla w, \phi_t +\bu\cdot\nabla \phi)_{\Omega}
+(\nabla_{\xi} w, \kappa_c \nabla_{\xi} \phi)_{\Omega}=0  \qquad \forall w \in \mathcal{V}^h
\end{align}
\end{minipage}}
\end{center}
where $\mathcal{V}^h$ is an appropriate discrete space, $\tau$ is the stabilisation parameter computed as
\begin{align}
\tau = \left (\Delta t ^2 + \bu\cdot \B{G} \bu \right ) ^{-1/2},
\end{align}
and $\kappa_c$ is the discontinuity capturing parameter computed as
\begin{align}
\kappa_c  = C \left |\phi_t +\bu\cdot\nabla \phi \right |,
\end{align}
where $C$ is an $O(1)$ parameter. 
Note that the residual-based nature choice of the discontinuity capturing parameter leads to  a strongly  consistent method.
This discontinuity capturing is inspired on similar discontinuity capturing reported in 
\cite{Guermond2014198} 
and is dimensionally consistent with more traditional discontinuity capturing parameter  such as reported in for instance \cite{JoKn07,Guermond20114248}.

The equations are integrated in time with the  Crank-Nicolson method. 
The resulting nonlinear system is solved with PETSc \cite{petsc-web-page,petsc-user-ref}.

\subsection{Volume conservation}\label{sec:vol}

Often, for instance when dealing with incompressible fluids,  it is  known {\it a priori}  
that the volume of each subdomain, given by
\begin{align}
V_0 =& \int_{\Omega} 1- H(\hat{\phi}) d\Omega, \nonumber \\
V_1 =& \int_{\Omega}  H(\hat{\phi}) d\Omega,
\end{align}
is conserved.
This mass conservation is not guaranteed in the discrete setting.
Due to accumulation of volume errors one constituent can completely disappear in extreme situations. 
To avoid this we add a global constant to our convected level-set field 
$\phi$ which restores global conservation. 
This global constant follows from solving the following non-linear equation:
\begin{align}
 \int_{\Omega}  H(\hat{\phi}(\phi^{n+1}+ \phi'^{n+1})) d\Omega = \int_{\Omega}  H(\hat{\phi}(\phi^{n}+. \phi'^{n})) d\Omega 
\end{align}
A simple Newton-Raphson procedure is employed to find this solution. 
This routine is quite fast as this is a scalar equation.

For convenience we introduce the mass correction in the time derivative of the overall convection problem. 
The approximation for time derivative is given by
\begin{align}
\phi_t = \frac{\phi^{n+1} - (\phi^{n} +\phi'^n)}{\Delta t}.
\end{align}
In this way $\phi'^{n+1}$ as an instantaneous correction necessary for maintaining mass conservation.


\subsection{Isogeometric analysis}

Isogeometric analysis, introduced  in \cite{HuCoBa04}, is the idea 
to directly employ shape functions from CAD  in simulations.
Since its introduction a large body of research is directed towards this general concept, 
see for instance \cite{CoHuBa09} and references there in.
One of the more popular choices for CAD shape functions are the so-called NURBS \cite{PiegTil97}.
In this work  NURBS will be employed as well.

Approximation properties similar to those of  standard finite elements can be proved\cite{BBCHS06} .
Furthermore, NURBS show superior convergence properties in terms 
of accuracy per degree of freedom \cite{EvBaBaHu08}.
A large body of work exist confirming 
the benefits of  adopting NURBS shapefunctions -- or  isogeometric analysis in general.


\section{Numerical experiments}
\label{sec:NumExp}

\subsection {Distortion test}\label{sec:dist}
To assess the four alternative approaches for redistancing we apply these 
to distorted meshes. The meshsize 
is uniformly reduced toward both the $x$ and $y$ axes. The underlying level-set functions, given by,
\begin{align}
\phi  = x-y
\end{align}
result in an interface from the bottom left to the top right. 

\begin{figure}[!h]
\begin{center}
\subfloat[][{\it Direct redistancing}]
{
\includegraphics[width=3.9cm]{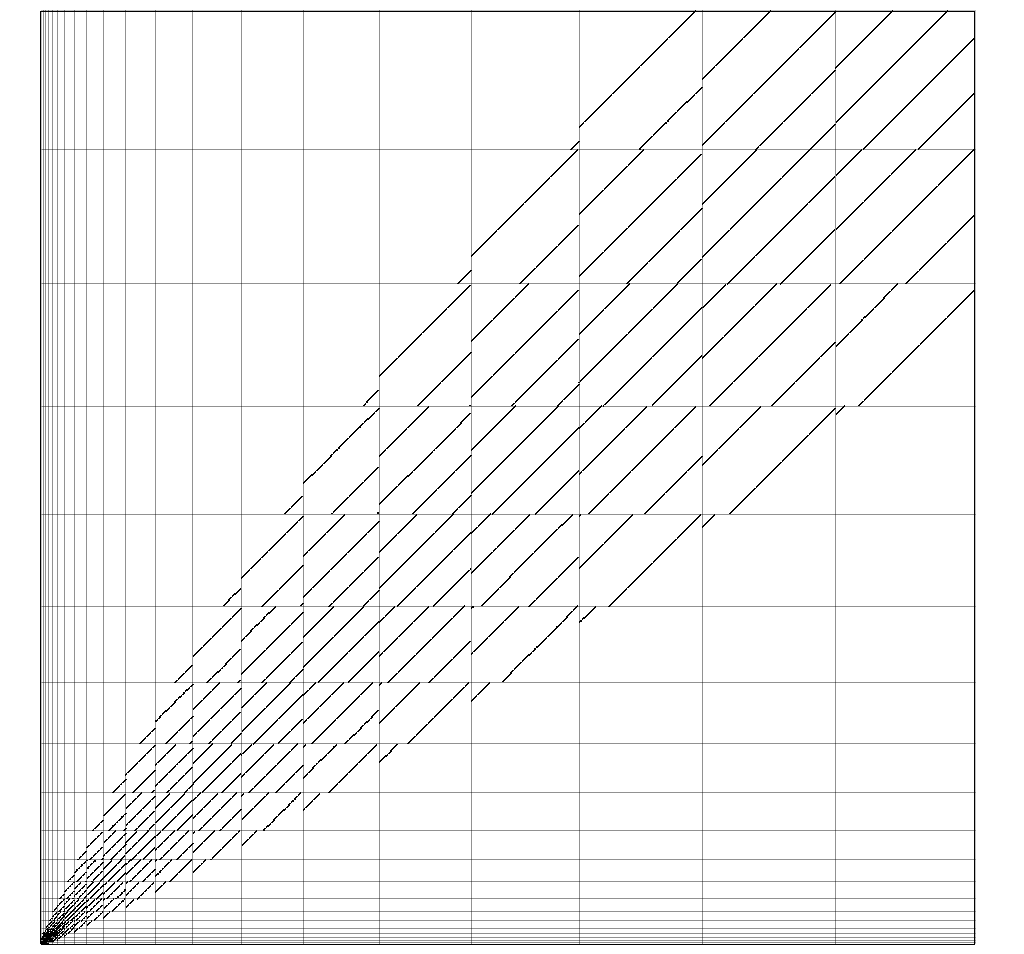}
}
\subfloat[][{\it Projected redistancing}]
{
\includegraphics[width=3.9cm]{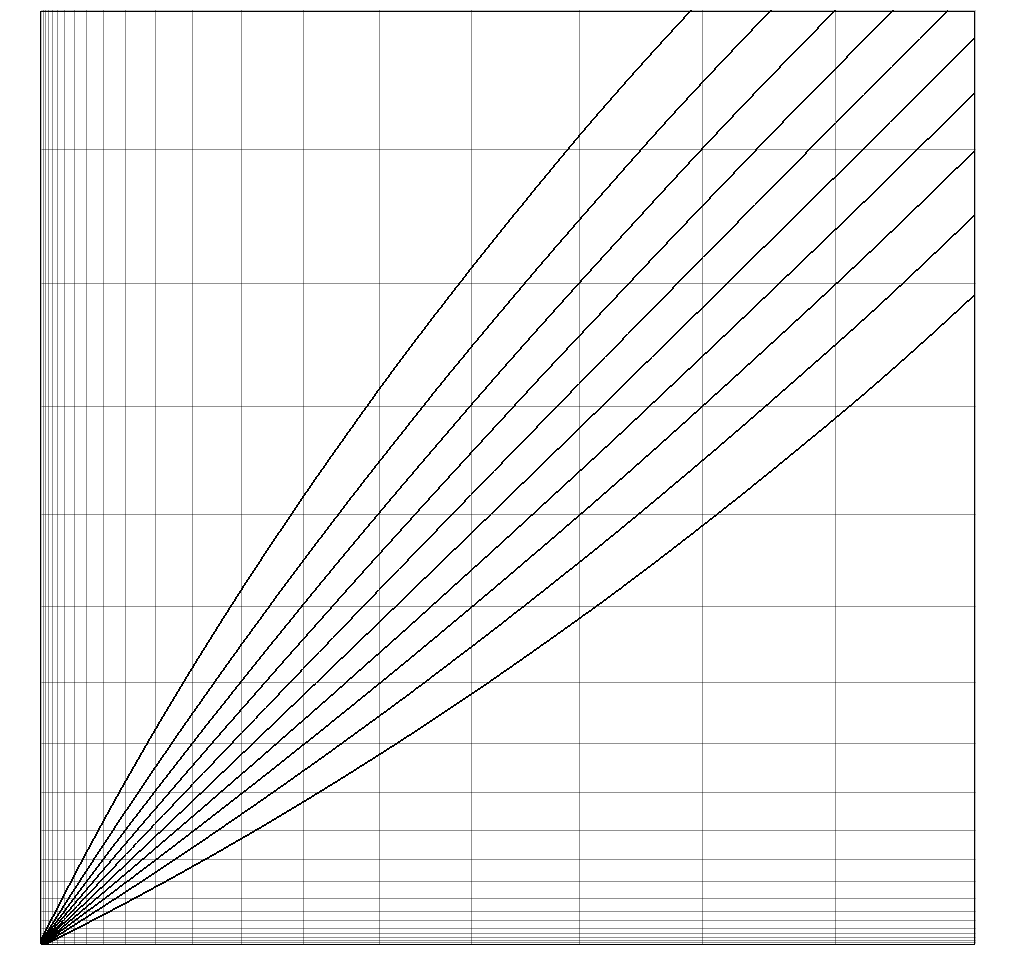}
}
\subfloat[][{\it Projected scaling}]
{
\includegraphics[width=3.9cm]{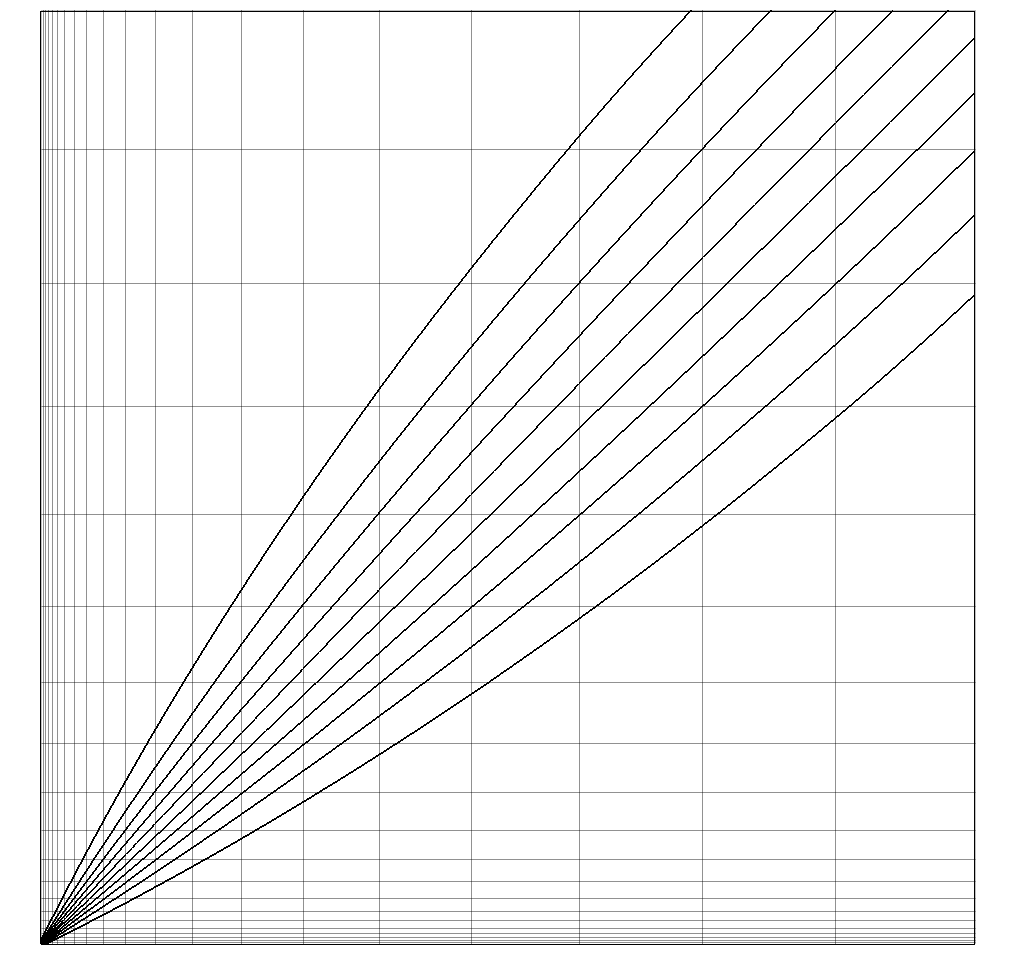}
}
\subfloat[][{\it Projected inverse scaling}]
{
\includegraphics[width=3.9cm]{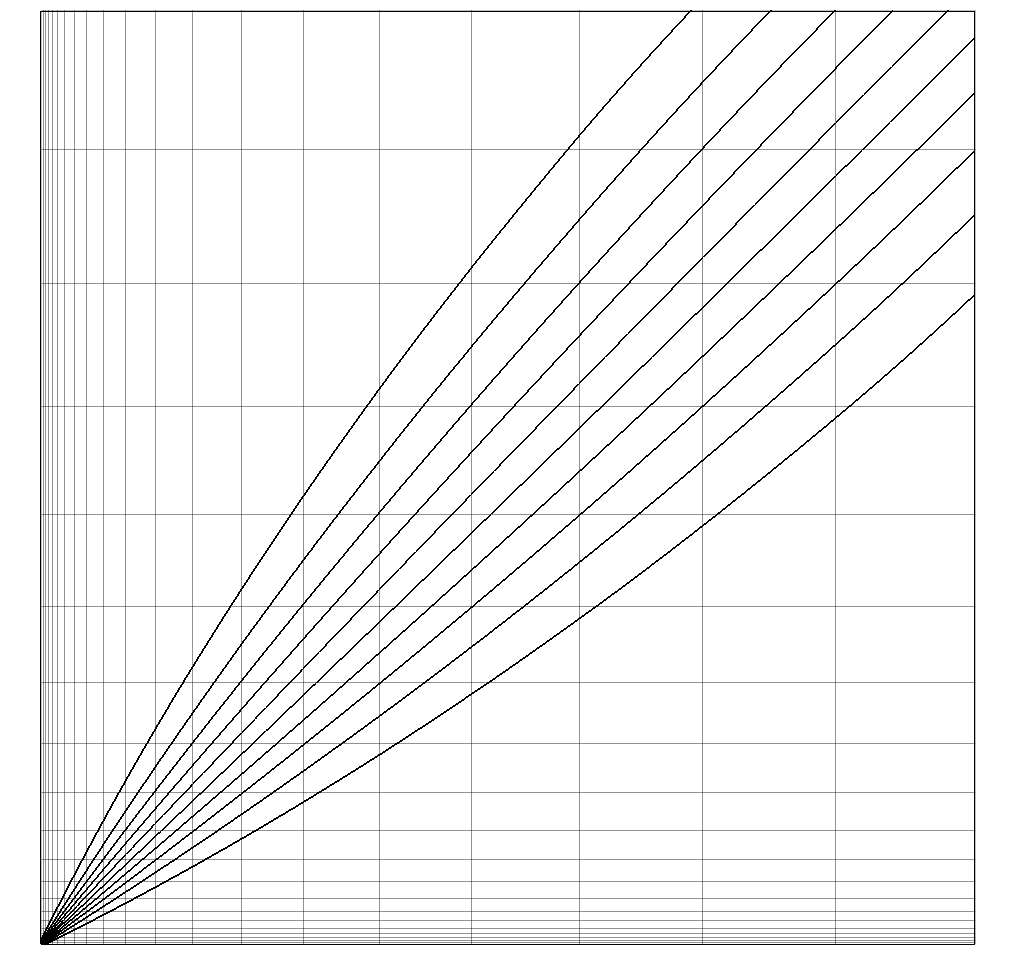}
}
\caption{Contour lines for the different redistancing approaches, 
all results are on a $C_1-Q_2$ mesh without any smoothing.} 
\label{fig:dist_k0}
\end{center}
\end{figure}

\begin{figure}[!h]
\begin{center}
\subfloat[][{\it Projected redistancing}]
{
\includegraphics[width=3.9cm]{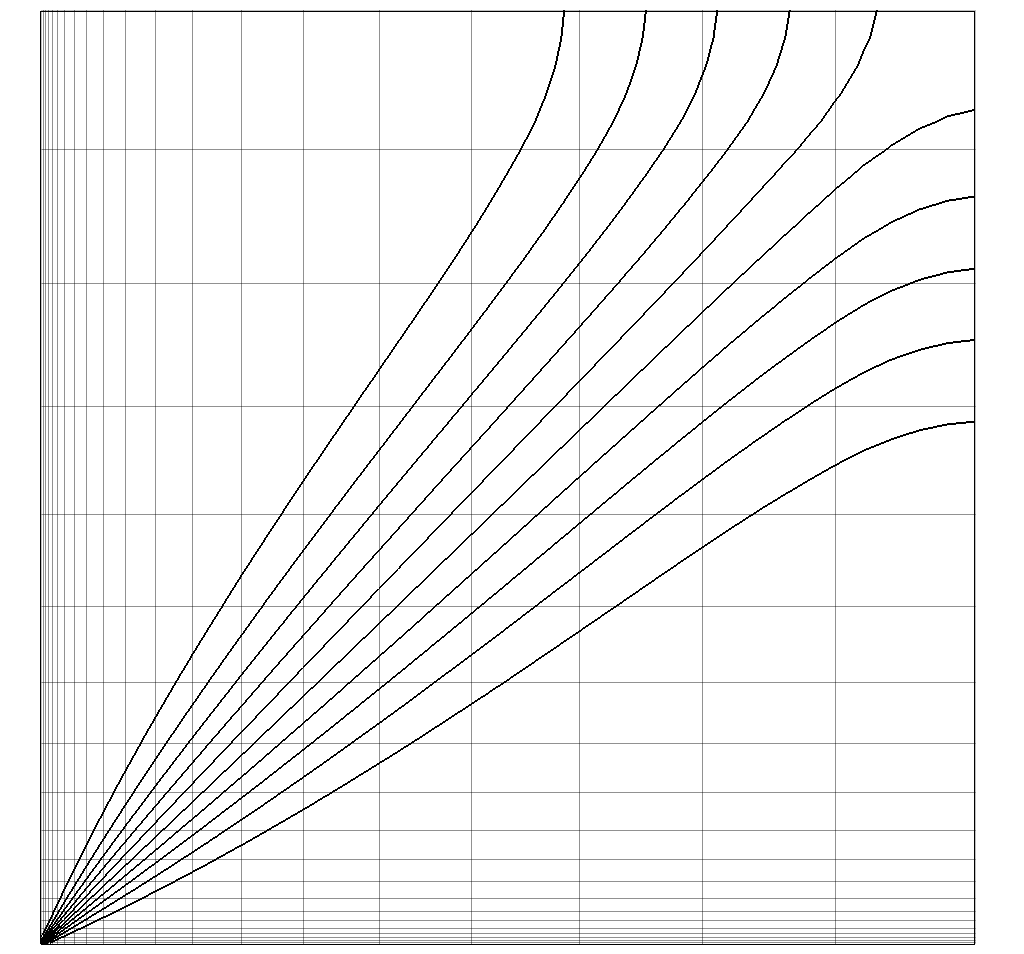}
}
\subfloat[][{\it Projected scaling}]
{
\includegraphics[width=3.9cm]{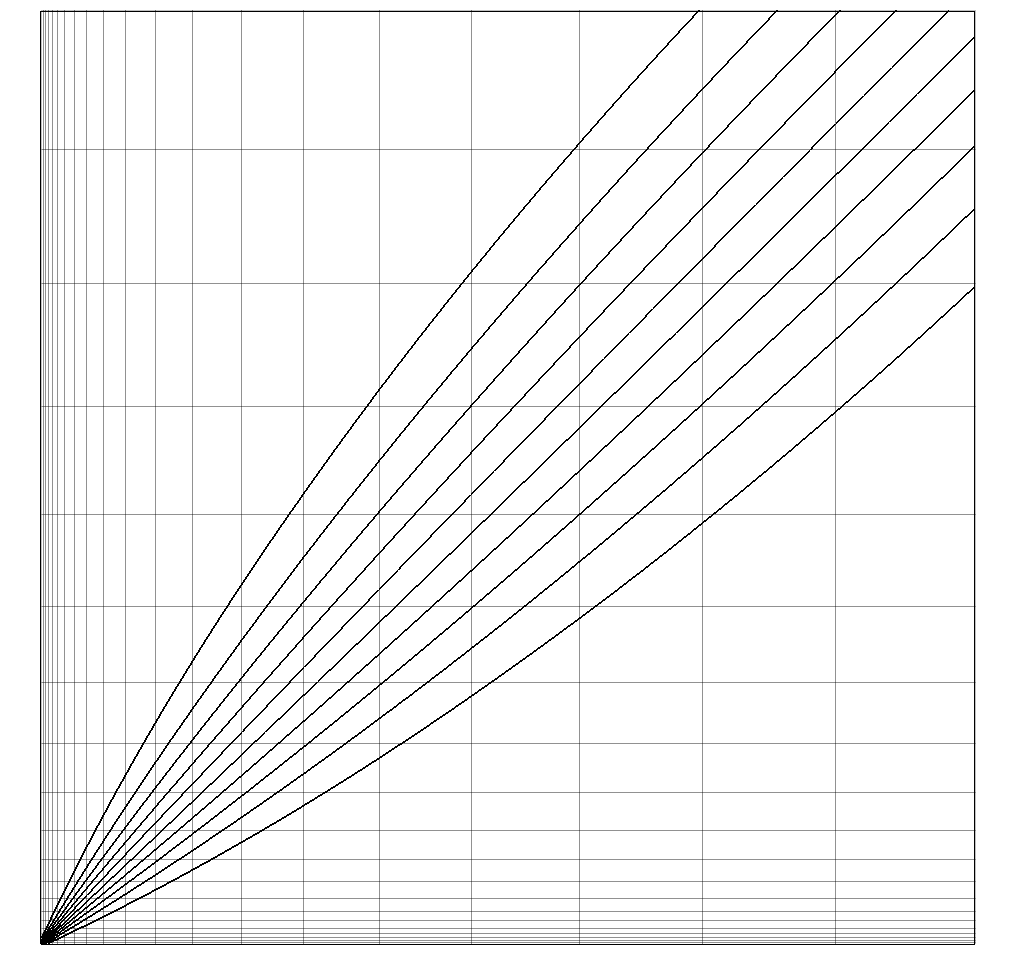}
}
\subfloat[][{\it Projected inverse scaling}]
{
\includegraphics[width=3.9cm]{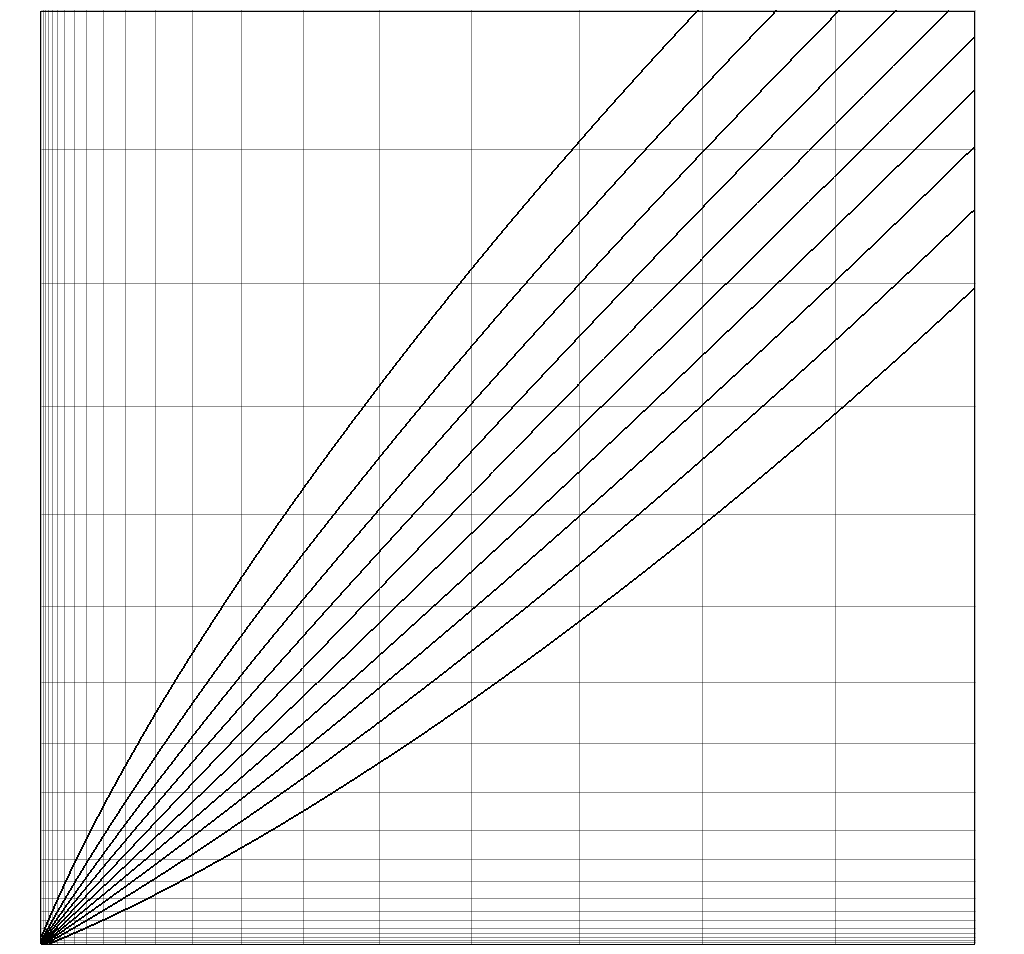}
}
\caption{Contour lines for the different projection approaches, all results are on a $C_1-Q_2$ mesh with smoothing $\kappa_d=1$.}
\label{fig:dist_k1}
\end{center}
\end{figure}

\begin{figure}[!h]
\begin{center}
\subfloat[][]
{
\includegraphics[width=3.9cm]{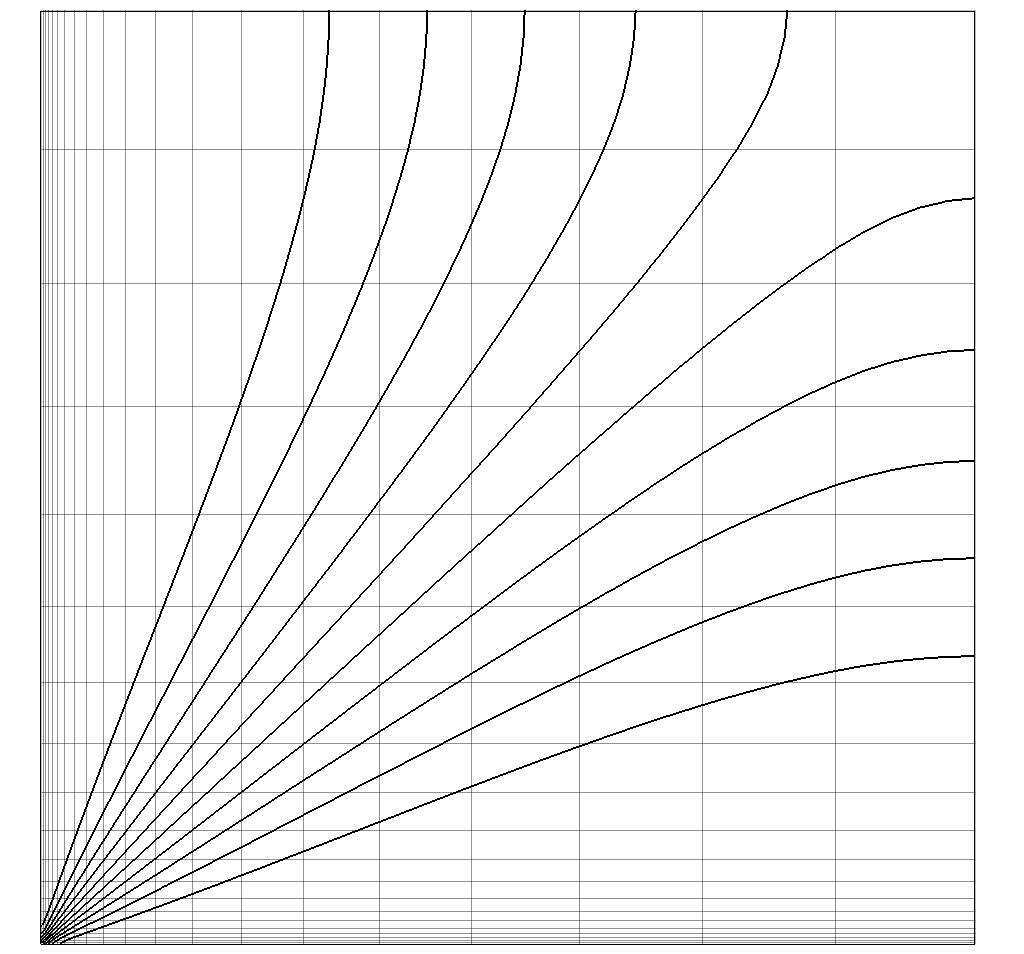}
}
\subfloat[][]
{
\includegraphics[width=3.9cm]{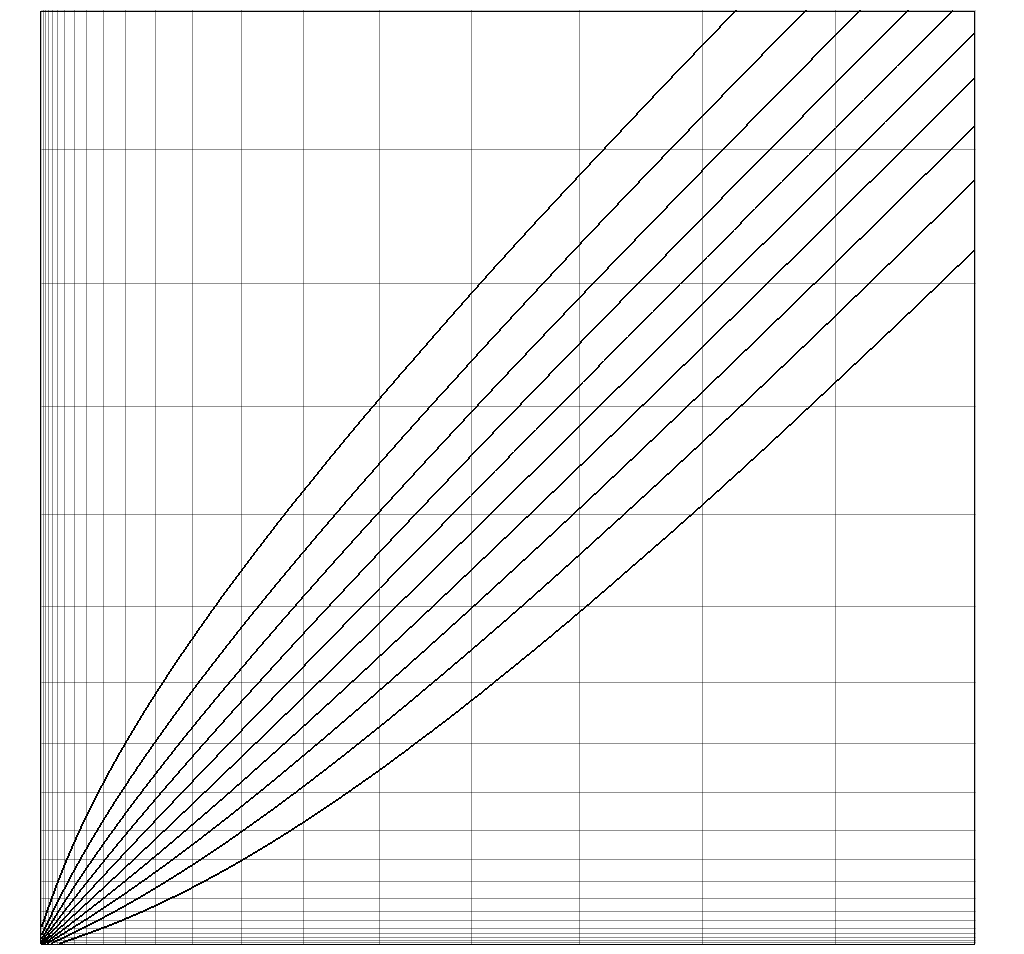}
}
\subfloat[][]
{
\includegraphics[width=3.9cm]{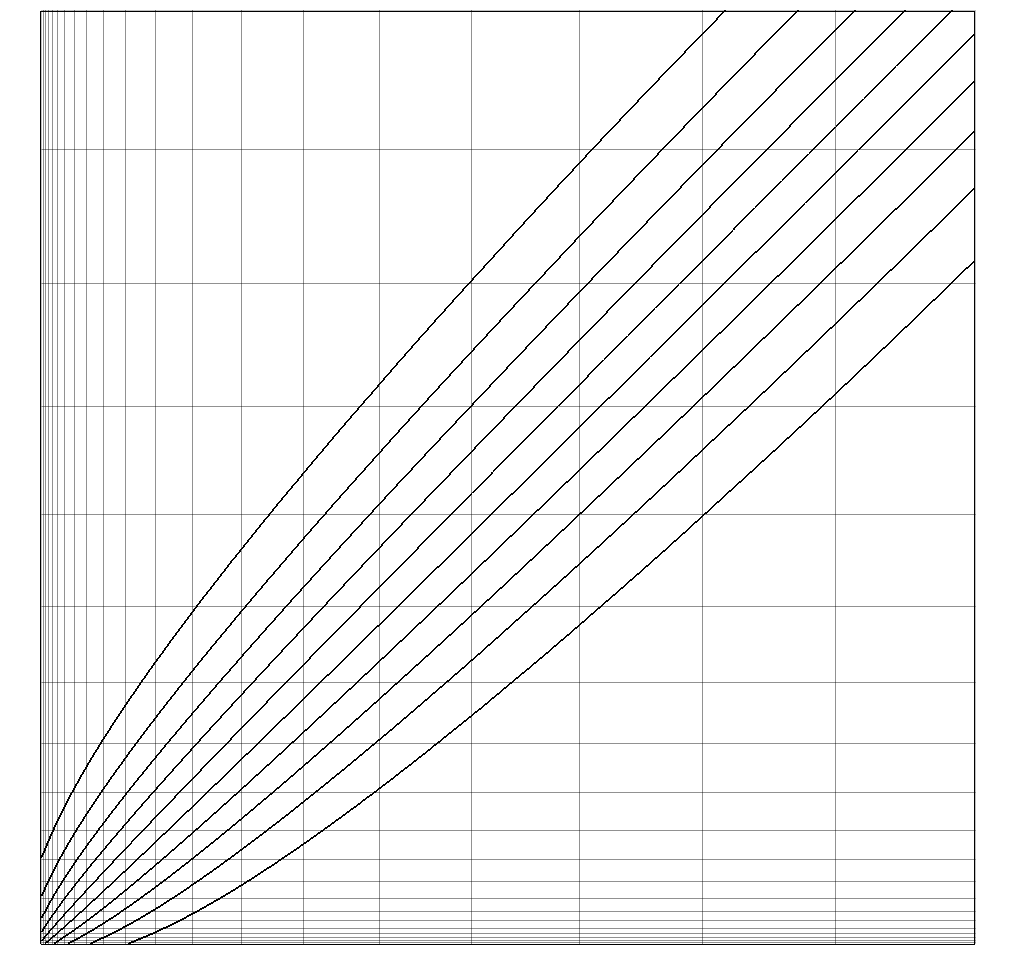}
}
\caption{Contour lines for the scaling approaches, all results are on a $C_1-Q_2$ mesh with
high smoothing $\kappa_d=10$.}
\label{fig:dist_k10}
\end{center}
\end{figure}

\begin{figure}[!h]
\begin{center}
\subfloat[][{\it Projected scaling}]
{
\includegraphics[width=5cm]{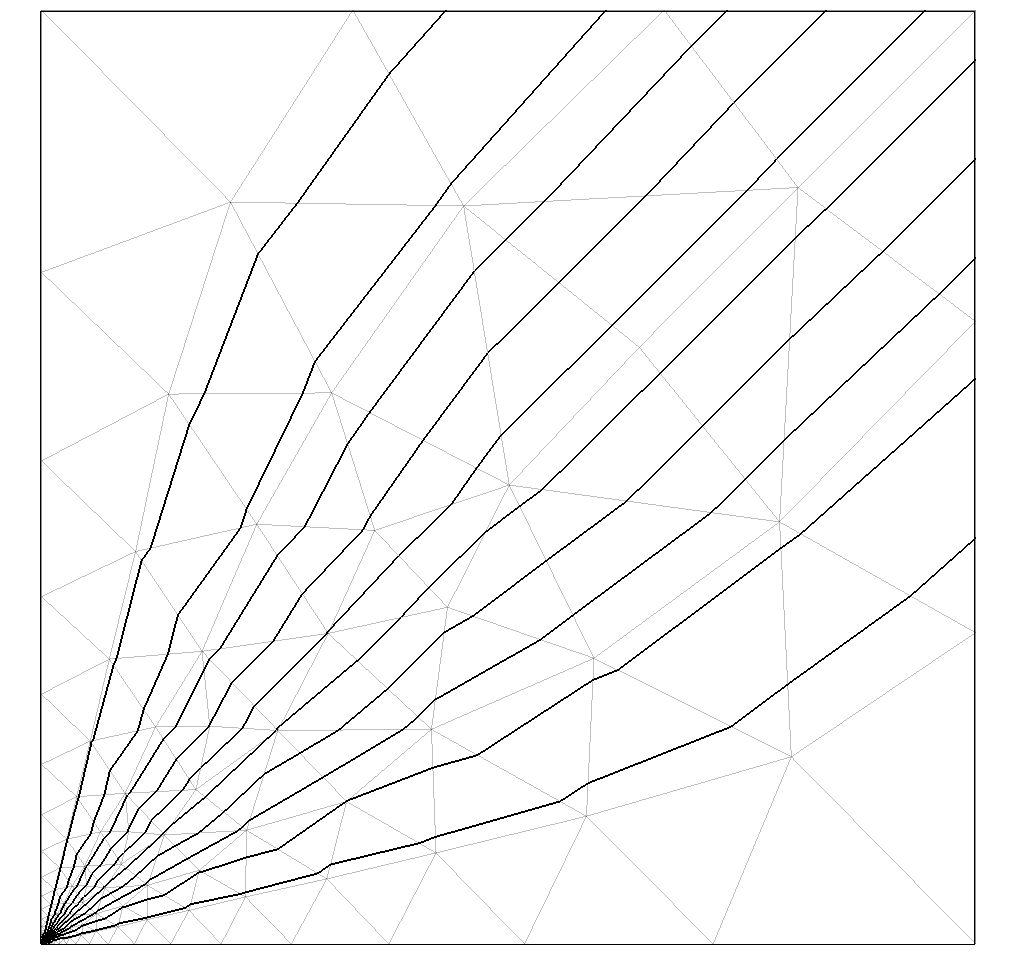}
}
\subfloat[][{\it Projected inverse scaling}]
{
\includegraphics[width=5cm]{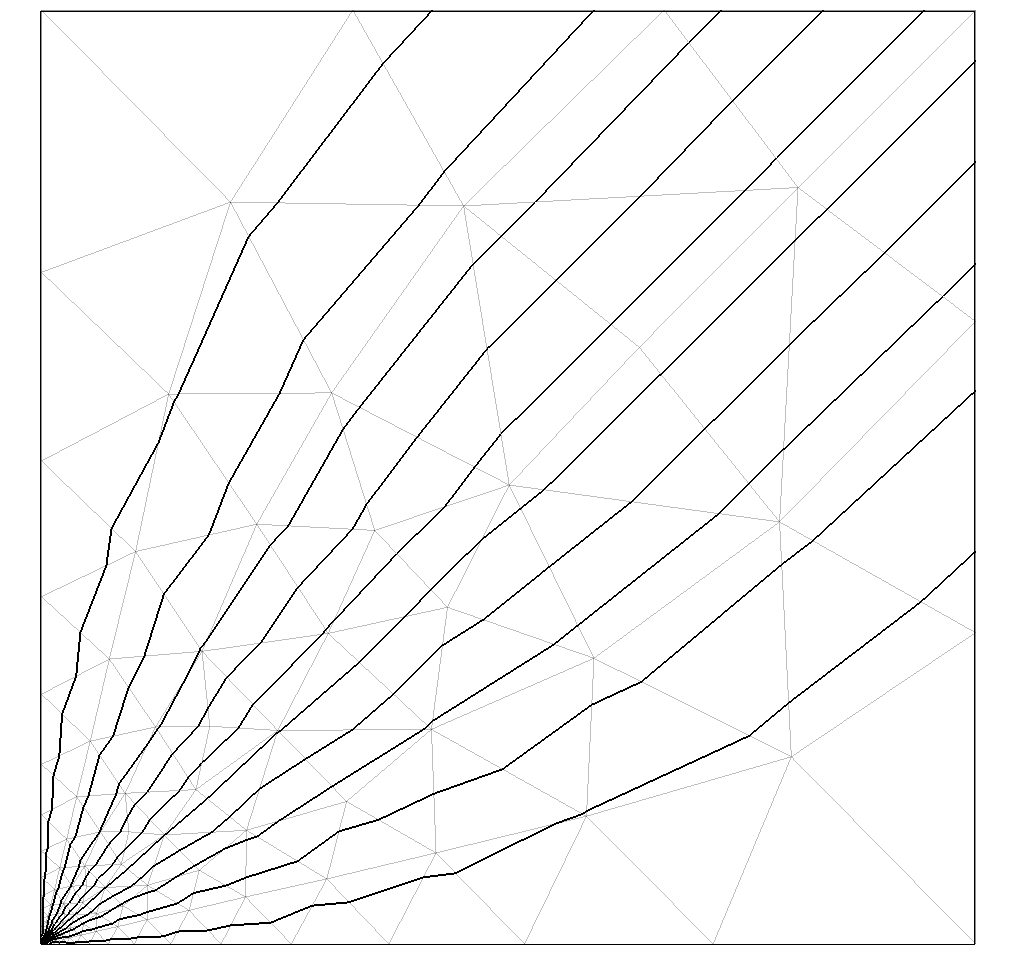}
}
\caption{Contour lines for the different redistancing approaches on a triangular mesh with smoothing $\kappa_d=1$.}
\label{fig:dist_tria}
\end{center}
\end{figure}

The contour lines of the Heaviside due to the redistancing in 
terms of the mesh length is given in figures \ref{fig:dist_k0},
\ref{fig:dist_k1},\ref{fig:dist_k10} and \ref{fig:dist_tria}. 
To exaggerate the results we have selected a relatively wide interface width, namely $\alpha=3$.

In figure \ref{fig:dist_k0} the contours for all four methods on a $C_1-Q_2$ mesh are presented. 
From this it can be clearly seen that the {\it direct redistancing} results in a discontinuous Heaviside, 
illustrated by the contour lines that terminate at element edges.  This shows that even on $C_1$ 
meshes the direct scaling does not necessarily provide a smooth Heaviside.  For the other three 
approaches the contour lines are nicely continuous and there is hardly 
any distinction between the three projection based methods.
On a $C_0-Q_1$ mesh the result for all four approaches will look analogous. 

In figure \ref{fig:dist_k1} the effect of moderate smoothing is presented. 
{\it Projected redistancing} clearly suffers from this inclusion of smoothing,
while the two scaled projection methods are largely unaffected. Although 
for {\it projected inverse scaling} the contour lines converge less at 
the origin due to the smoothing. With increased smoothing  {\it projected inverse scaling} 
suffers more from this effect, as shown in figure \ref{fig:dist_k10}.

In figure \ref{fig:dist_tria} the  two scaled projection methods 
on a triangular mesh are presented.
The triangular meshes are generated with gmsh \cite{gmsh}
 In this case the irregularity of the triangular meshes highly favours the 
inclusion of some modest smoothing.

\subsection{Problem description: Vortex in a box}

To assess the efficacy of the proposed redistancing approaches we  
apply all methods on the vortex in a box problem. We focus only 
on the three projection-based alternatives, as the direct approach fails to deliver smooth Heaviside functions.

The problem setup originates from \cite{HighResLeveque} while the parameters 
are taken from \cite{Rider_Kothe_98}. The computational domain is a unit square
and $\phi$ is initialized such that its zero level set is a circular disc of radius $0.15$ 
centered at $(0.50, 0.75)$. The time-dependent velocity  $\bu = (u,v)$, is given by
\begin{align}
u =& \cos \left (\frac{\pi t}{8} \right ) \sin (2\pi y) \sin(\pi x)^2,  \nonumber \\
 v =&  -\cos \left (\frac{\pi t}{8} \right ) \sin (2\pi x) \sin(\pi y)^2. 
 \label{eq:2Ddef}
\end{align}
This velocity field stretches the disc into a very thin spiral until $t = 4$, next the flow reverses and the spiral is deformed back into a perfect circle at $t = 8$. 

We solve this problem on different meshes. On the one hand we have a sequence of 
quadrilateral NURBS elements 
 with linear-$C_0$  or quadratic-$C_1$  tensor product based shape functions. 
Additionally, we solve the problem on triangular elements, which indicates the generality of the approach.

We select an $80\times80$ element NURBS meshes and an equivalent  triangular 
mesh with $16794$ linear elements and $8558$ nodes. 
For the triangular mesh the target line element count along an domain edge is 
identical to the that of the NURBS meshes. 
This mesh size is selected because it is sufficiently accurate to draw meaningful conclusions, 
while on the other hand it is still coarse enough to see the numerical errors and judge the
 behaviour of the methods with the naked eye.

\subsection{ Projected  redistancing}
\begin{figure}[!h]
\begin{center}
\subfloat[][$\kappa_d=0.1$]
{
\includegraphics[width=5cm]{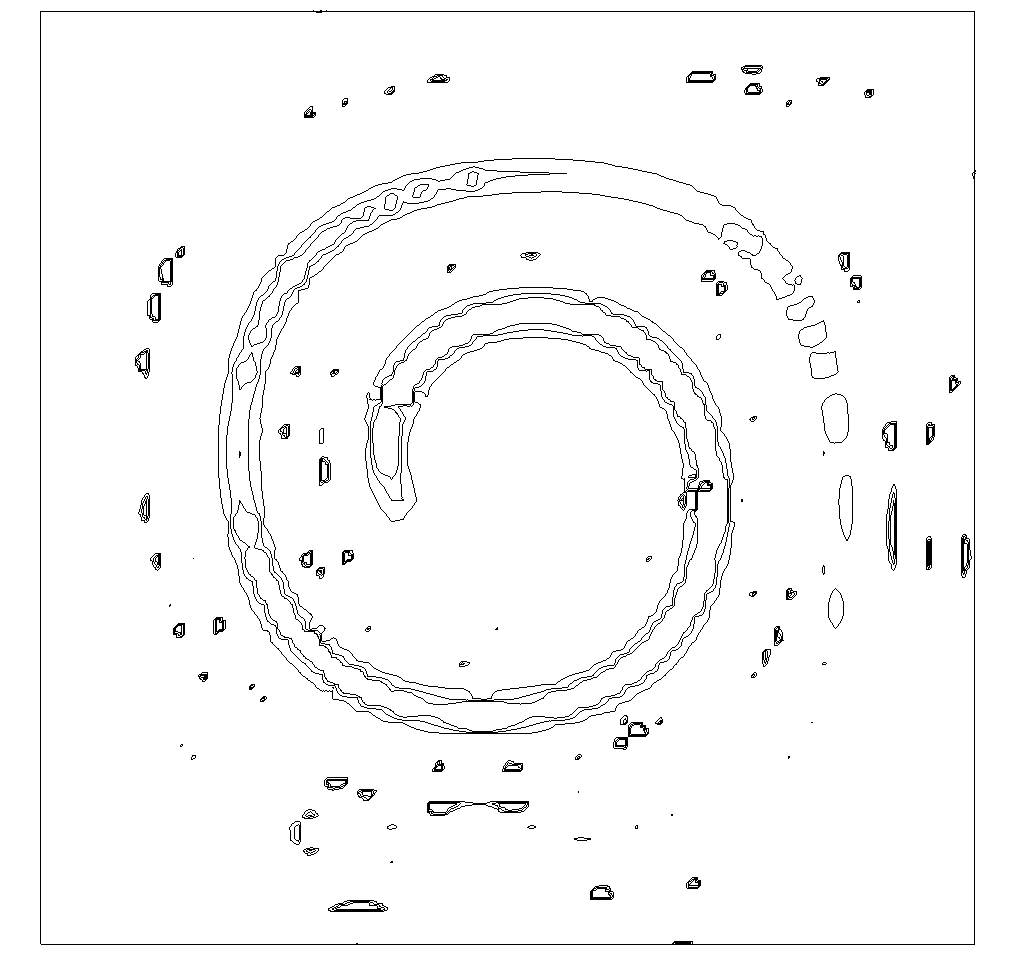}
}
\subfloat[][$\kappa_d=1$]
{
\includegraphics[width=5cm]{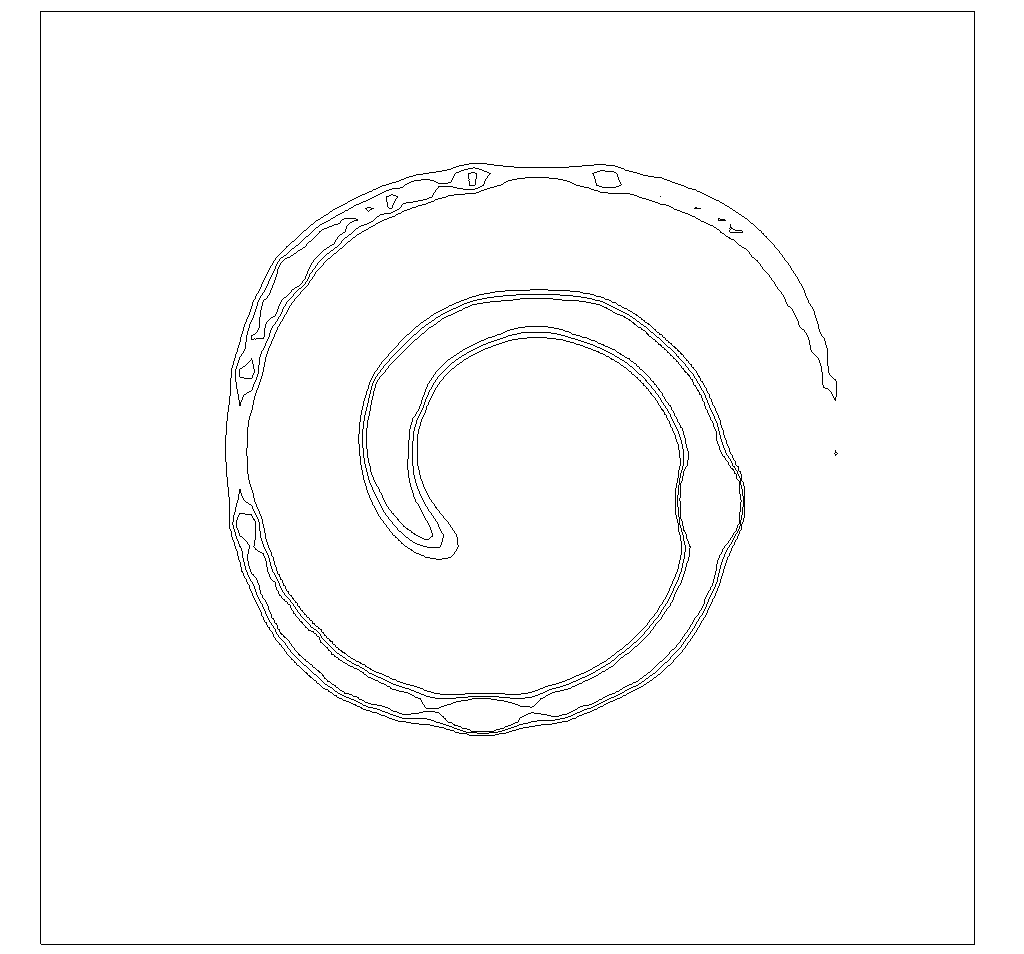}
}
\subfloat[][$\kappa_d=10$]
{
\includegraphics[width=5cm]{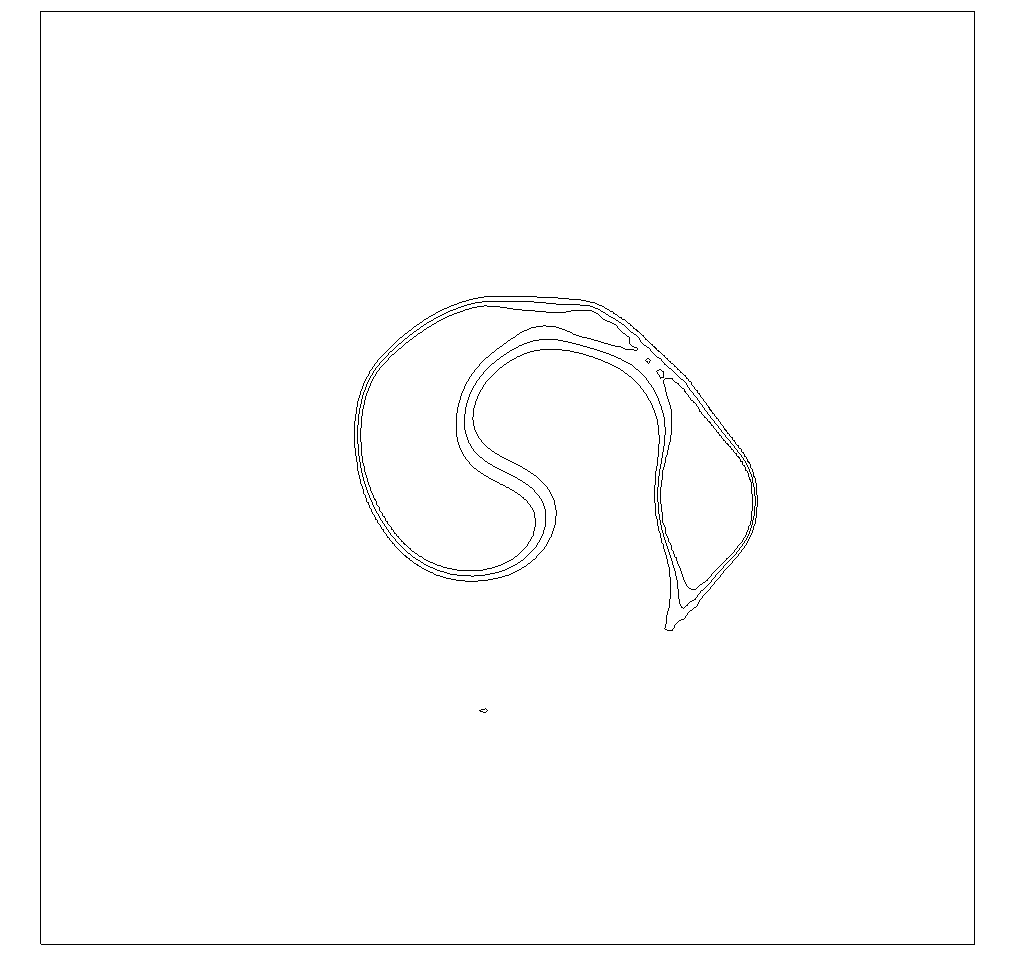}
}
\caption{Contour plots of the regularized Heaviside for {\it projected redistancing} on a $80 \times 80$ $C_0$ linear mesh.}
\label{fig:mode2_c0}
\end{center}
\end{figure}

%
Figure \ref{fig:mode2_c0} shows 
snapshots of the regularized Heaviside at maximum  distortion for the different 
smoothing parameters $\kappa$ for a  $80 \times 80$ $C_0$ linear mesh.

It can be seen there is a delicate balance 
between not enough and to much smoothing. When $\kappa_d$  is 
small there is not enough smoothing and the distance field shows oscillations, 
while for the  higher values of $\kappa_d$ the diffusion caused by the smoothing 
moves the interface more than desired.
For the other discretizations --  $P_1$  triangles  and $C_1-Q_2$ 
quadrilaterals  -- we get similar results, leading to the same conclusion. 

\begin{figure}[!h]
\begin{center}
\input{mode2/m2_k1_80p1_corr.tex}
\caption{Correction required for volume conservation for the $80\times80$ linear mesh with $\kappa_d=1$.}
\label{fig:mode2_corr}
\end{center}
\end{figure}

For all simulations the volume is conserved up to machine precision, even when 
the solution is highly deformed due to the smoothing.
This demonstrates effectivity of the  volume correction method presented  in section \ref{sec:vol}.
The correction required for volume conservation is of the  order $10^{-3}$ per time step. 
A typical time trace of this correction is given in 
figure \ref{fig:mode2_corr}.

%
\subsection{Projected scaling}
\begin{figure}[!h]
\begin{center}
\subfloat[][$\kappa_d=0.1$]
{
\includegraphics[width=5cm]{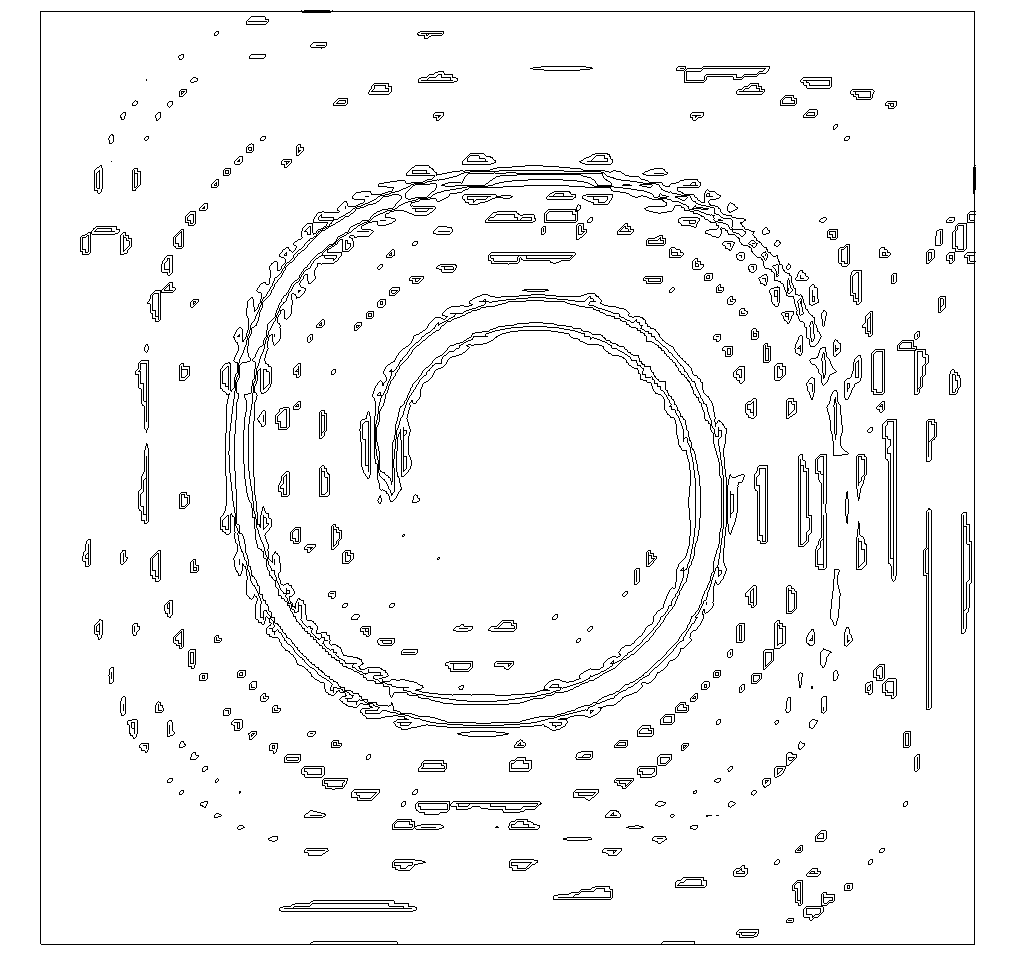}
}
\subfloat[][$\kappa_d=1$]
{
\includegraphics[width=5cm]{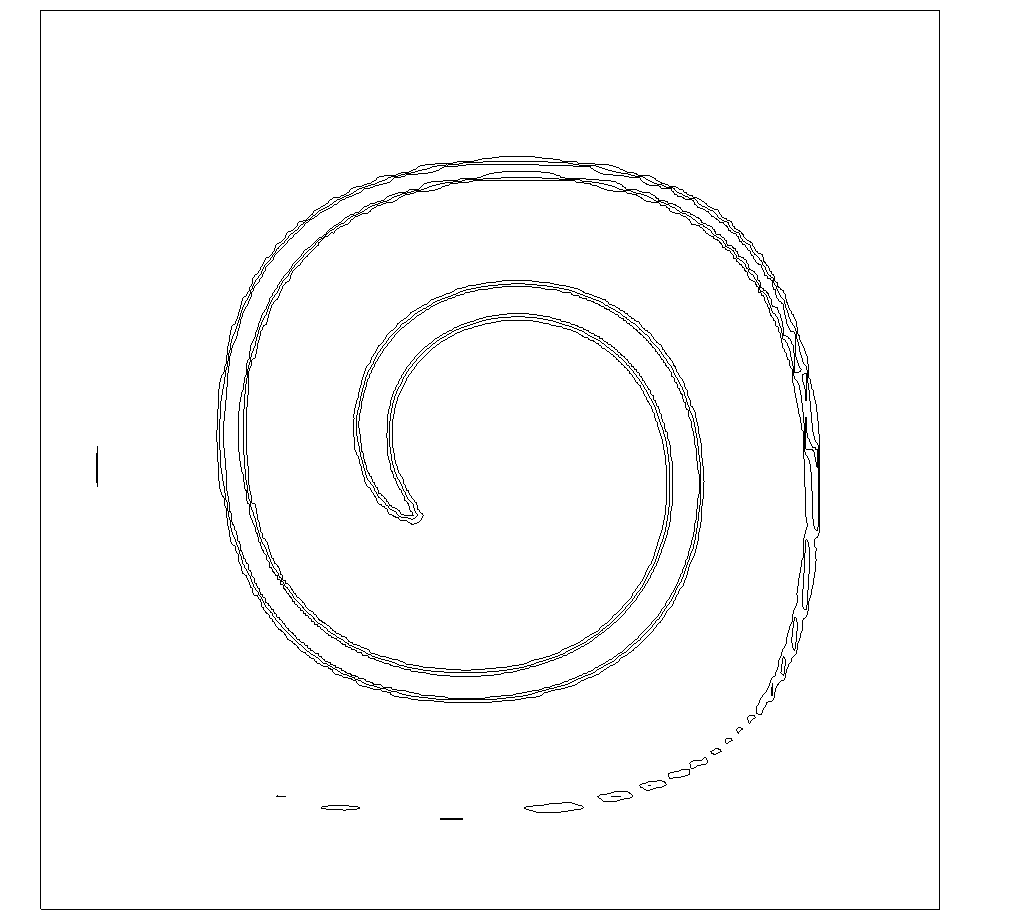}
}
\subfloat[][$\kappa_d=10$]
{
\includegraphics[width=5cm]{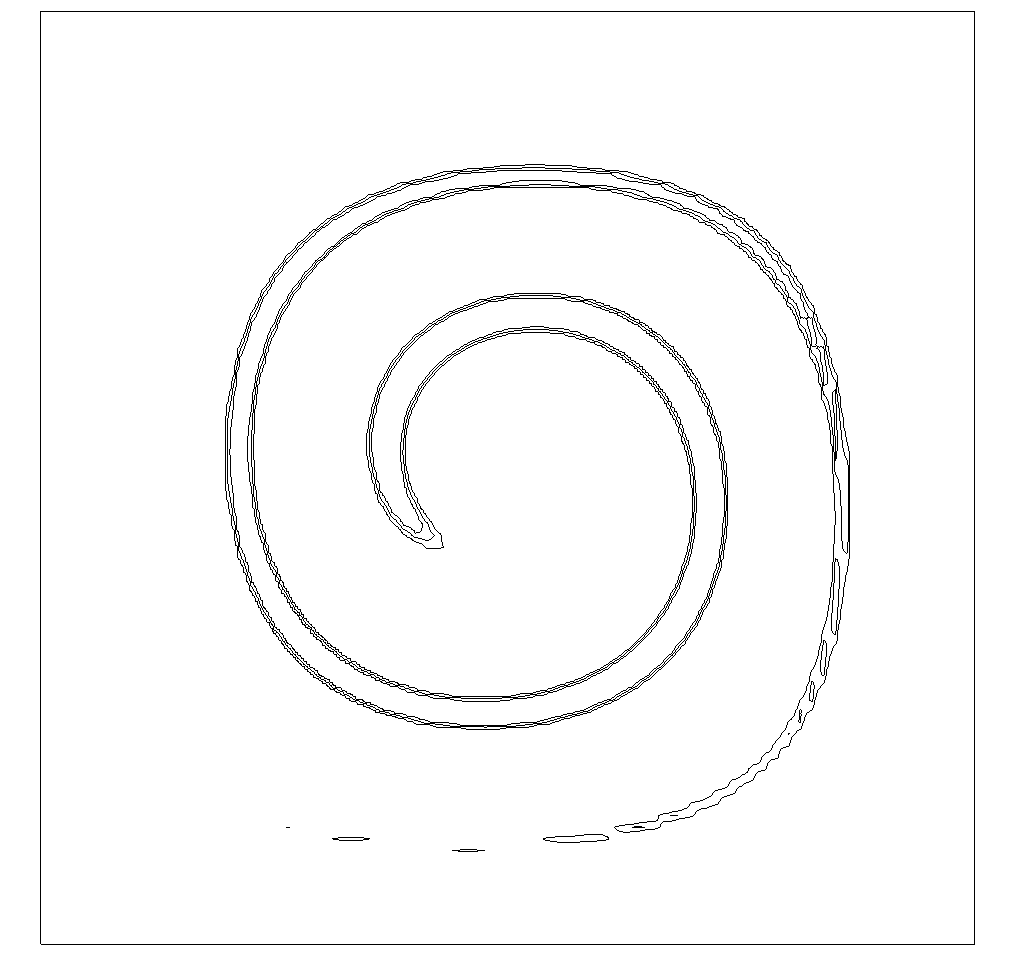}
}
\caption{Contour plots of the regularized Heaviside for {\it projected scaling} 
on a $80\times80$ $C_0$ linear mesh.}
\label{fig:mode3_c0}
\end{center}
\end{figure}

Figure \ref{fig:mode3_c0}  shows 
snapshots of the regularized Heaviside at maximum  distortion for the 
different smoothing parameters $\kappa_d$ for
a  $80 \times 80$ $C_0$ linear mesh.

It can be seen that for this approach there is no problem when a large smoothing parameter is used.
However, when the smoothing parameter is too small the solution becomes spotty 
and irregular similar as in the previous section.
This monotone behaviour with respect to the smoothing parameter makes this a viable approach.
For the other discretizations --  $P_1$  triangles  and $C_1-Q_2$ 
quadrilaterals  -- we get similar results, leading to the same conclusion. 

\begin{figure}[!h]
\begin{center}
\input{mode3/m3_k1_80p1_corr.tex}
\caption{Correction required for volume conservation for the $80\times80$ linear mesh with $\kappa_d=1$.}
\label{fig:mode3_corr}
\end{center}
\end{figure}

The interface conservative nature of the redistancing approach is also apparent in 
the time trace of the correction parameter required to conserve volume. A typical 
time trace is depicted in figure \ref{fig:mode3_corr}.  Comparing with the time 
trace of the previous section, in figure \ref{fig:mode2_corr},  the required 
correction is roughly an order of magnitude less.

%
\subsection{Projected inverse scaling} 
\begin{figure}[h]
\begin{center}
\subfloat[][$\kappa_d=0.01$]
{
\includegraphics[width=5cm]{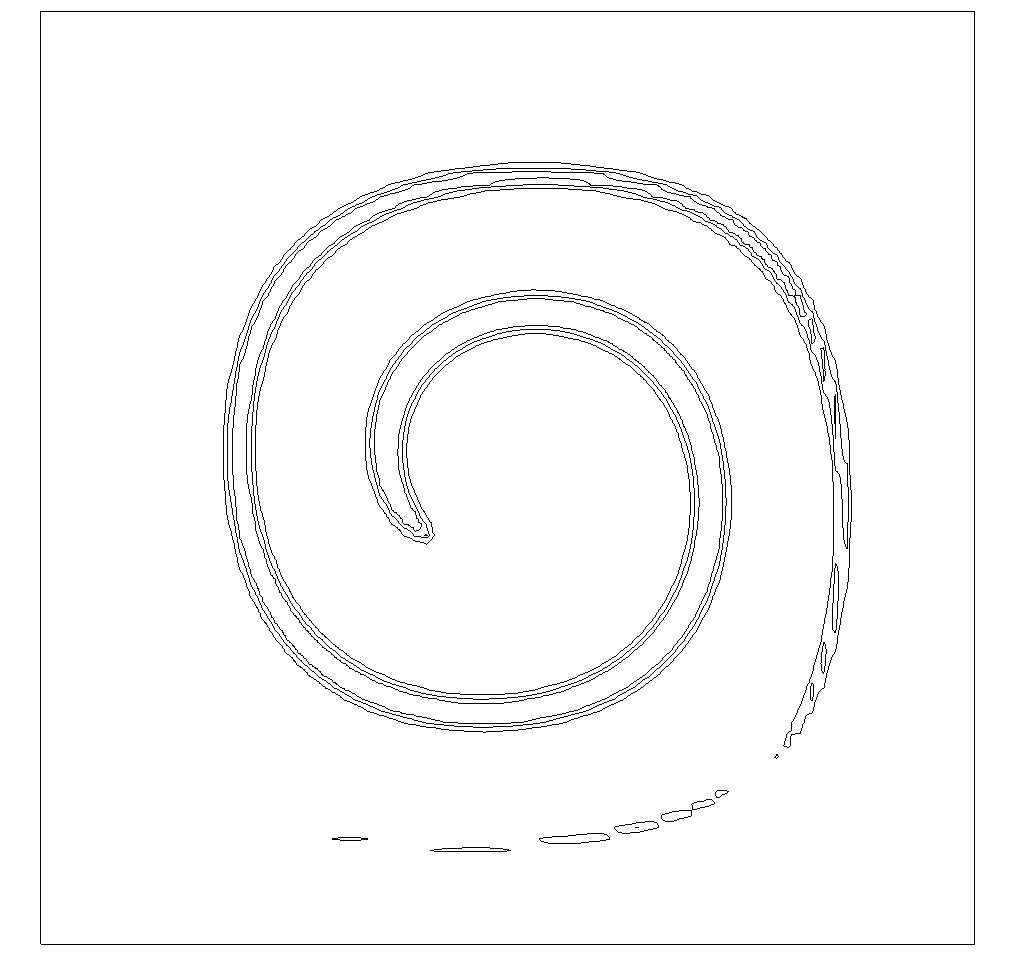}
}
\subfloat[][$\kappa_d=1$]
{
\includegraphics[width=5cm]{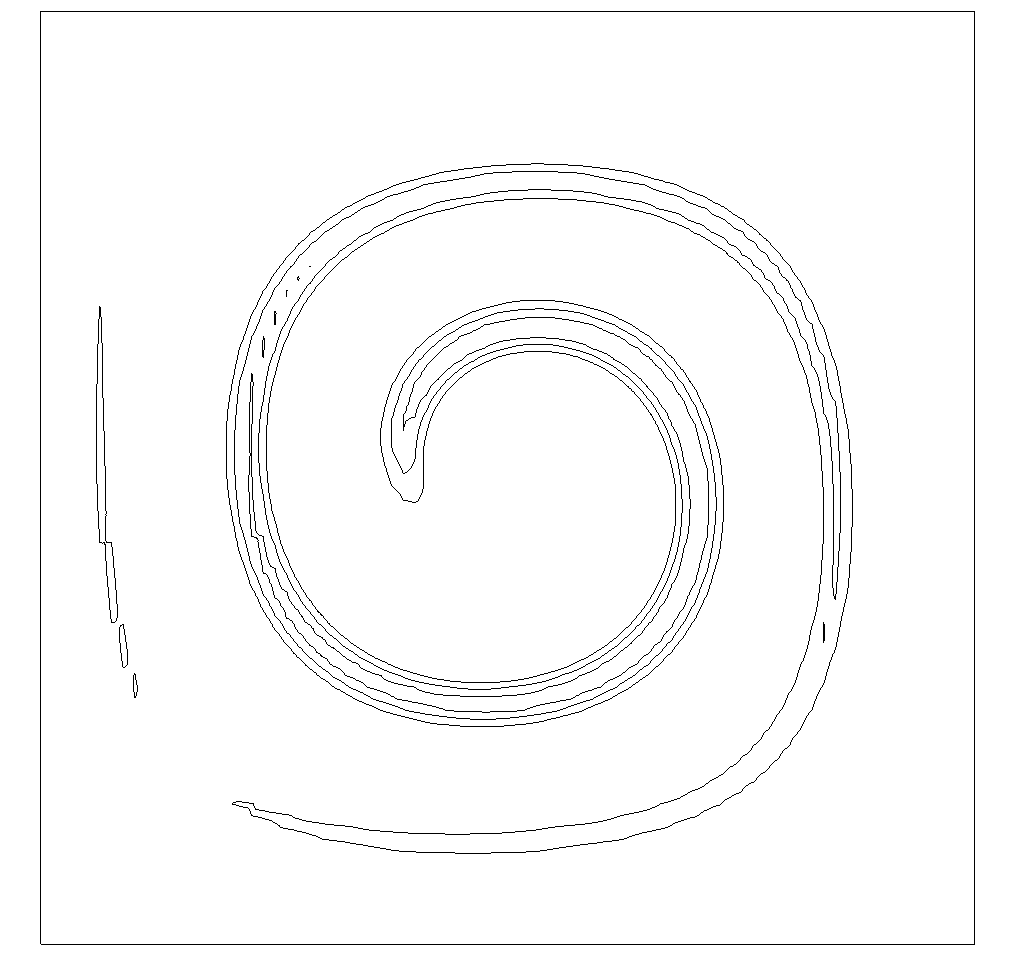}
}
\subfloat[][$\kappa_d=100$]
{
\includegraphics[width=5cm]{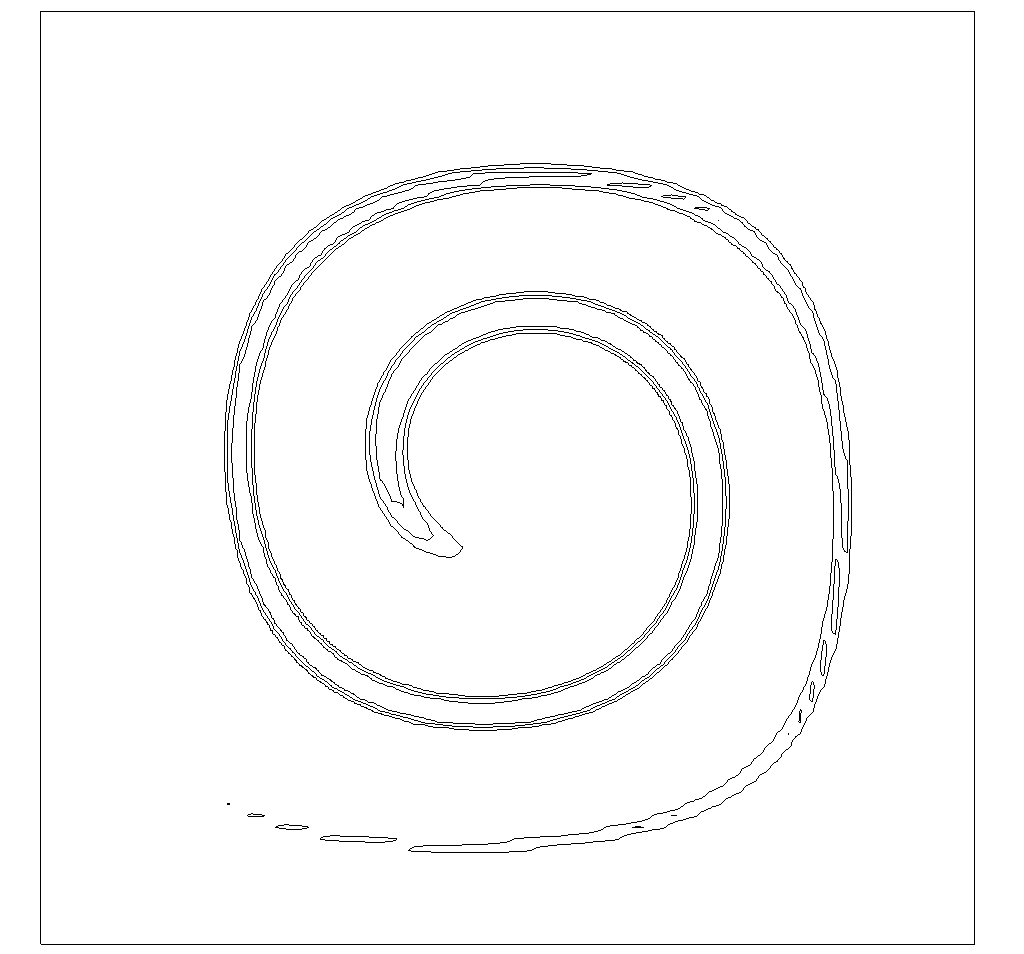} 
}
\caption{Contour plots of the regularized Heaviside for {\it projected inverse scaling} on a $80\times80$ $C_0$ linear mesh.}
\label{fig:mode4_c0}
\end{center}
\end{figure}

Figure \ref{fig:mode4_c0}  shows 
snapshots of the regularized Heaviside at maximum  distortion for the 
different smoothing parameters $\kappa$ for
a  $80\time80$ $C_0$ linear mesh. 
The smoothing does 
not have a significant influence on the  results. Both high and low smoothing limits  display acceptable answers.

\begin{figure}[!h]
\begin{center}
\subfloat[][$\kappa_d=0$]
{
\includegraphics[width=5cm]{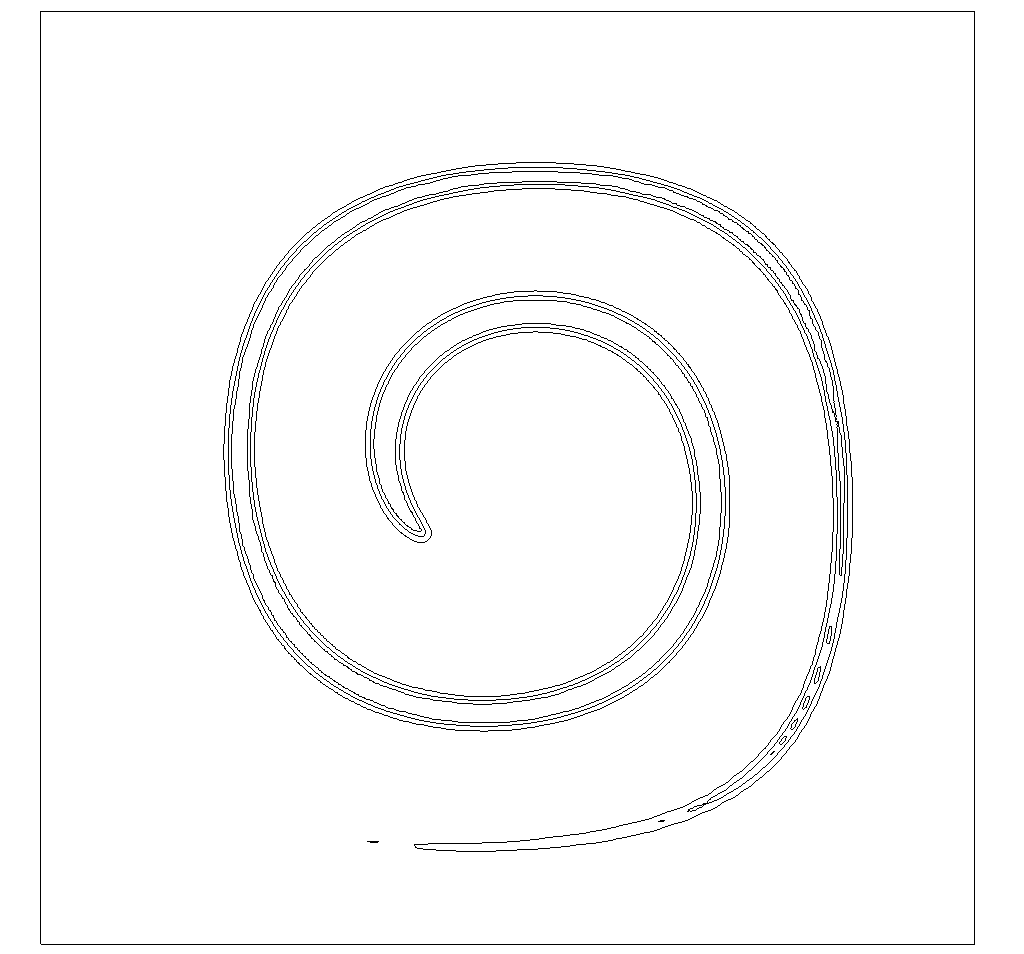}
}
\subfloat[][$\kappa_d=1$]
{
\includegraphics[width=5cm]{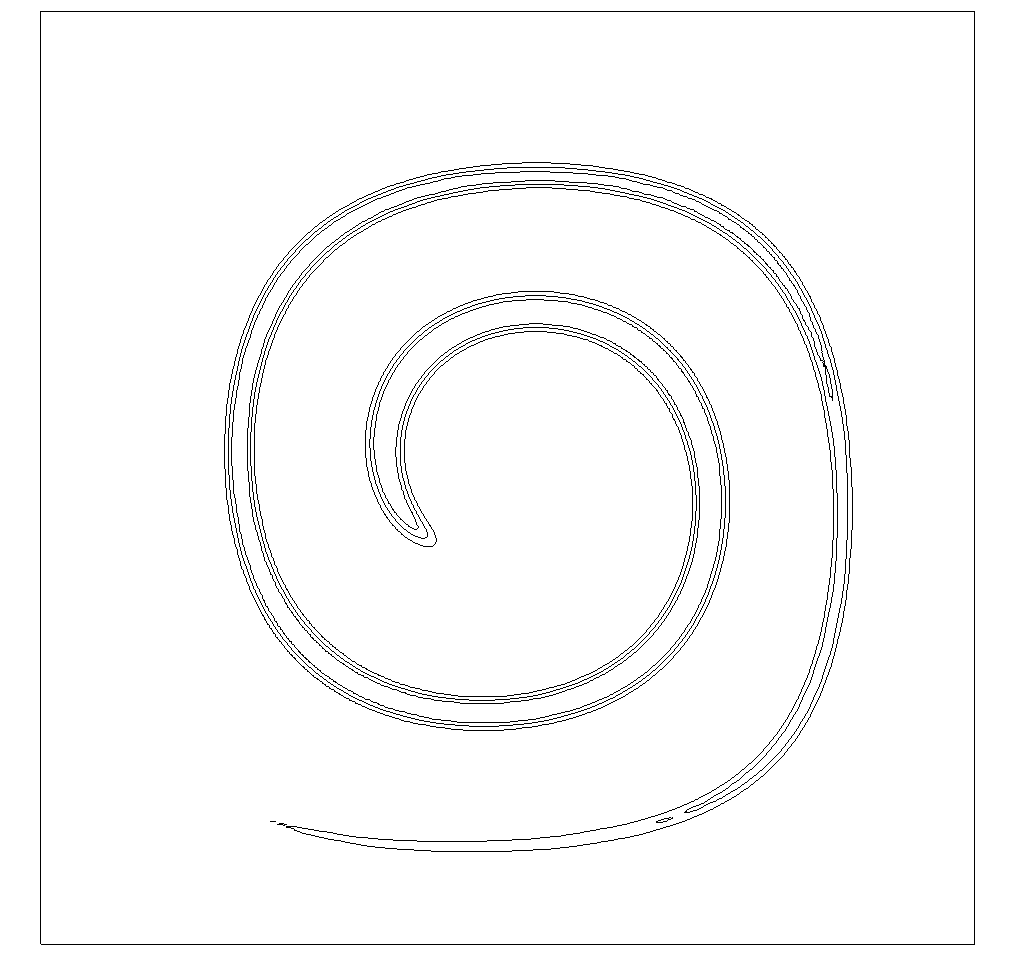}
}
\subfloat[][$\kappa_d=100$]
{
\includegraphics[width=5cm]{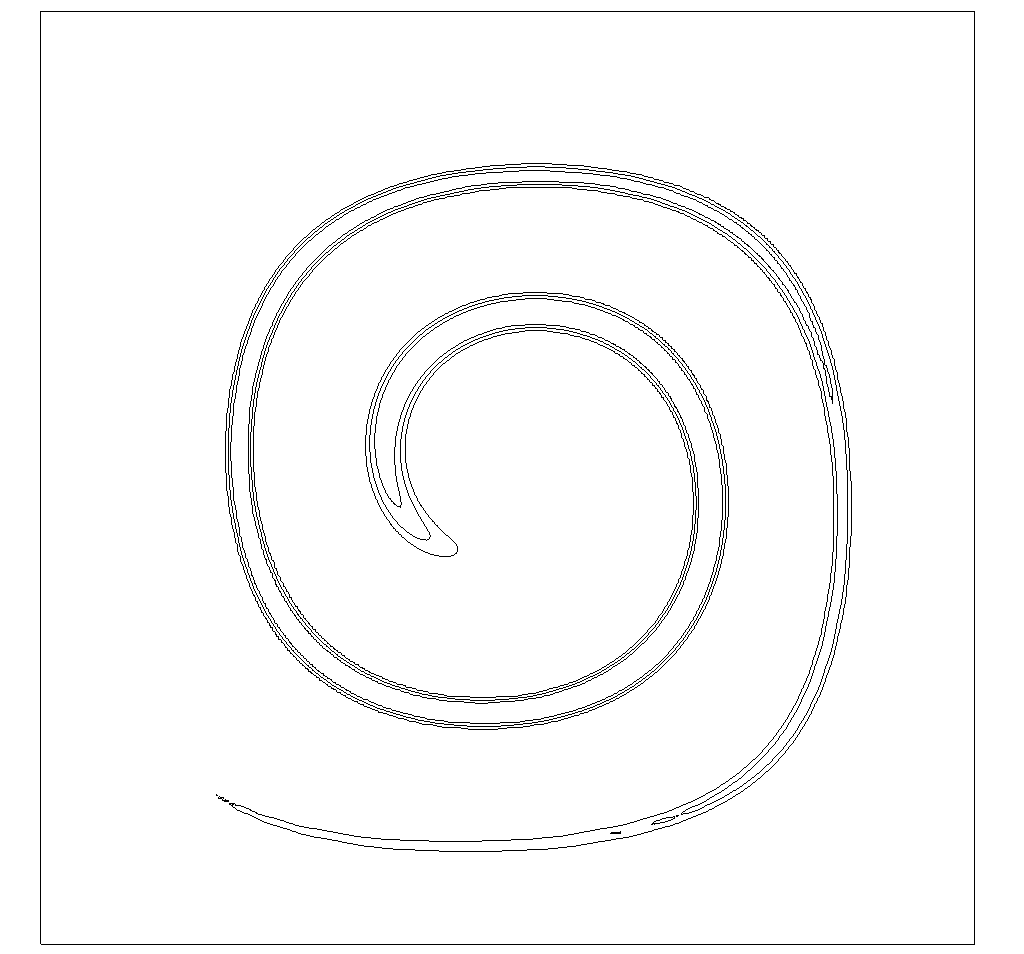}
}
\caption{Contour plots of the regularized Heaviside for {\it projected inverse scaling} on  a $80\times80$ $C_1$ quadratic mesh.}
\label{fig:mode4_c1}
\end{center}
\end{figure}

\begin{figure}[!h]
\begin{center}
\subfloat[][$\kappa_d=0.01$]
{
\includegraphics[width=5cm]{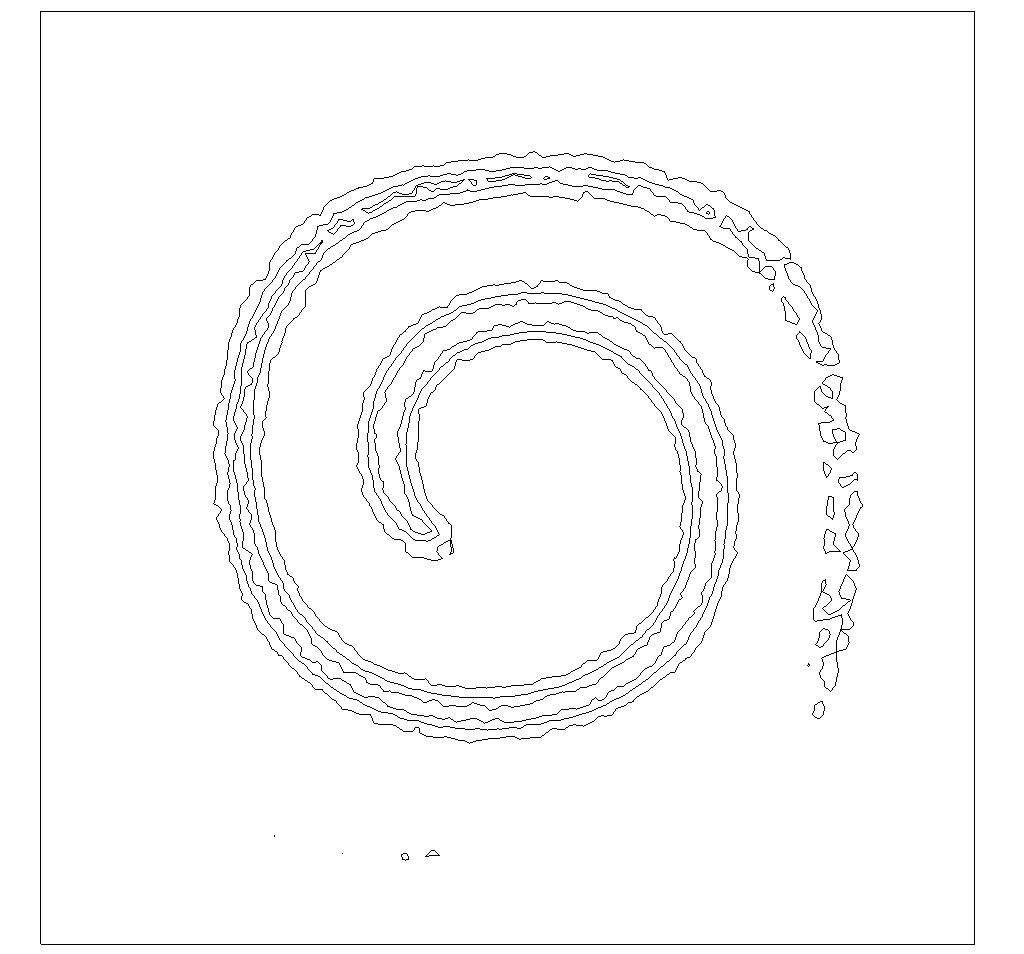}
}
\subfloat[][$\kappa_d=1$]
{
\includegraphics[width=5cm]{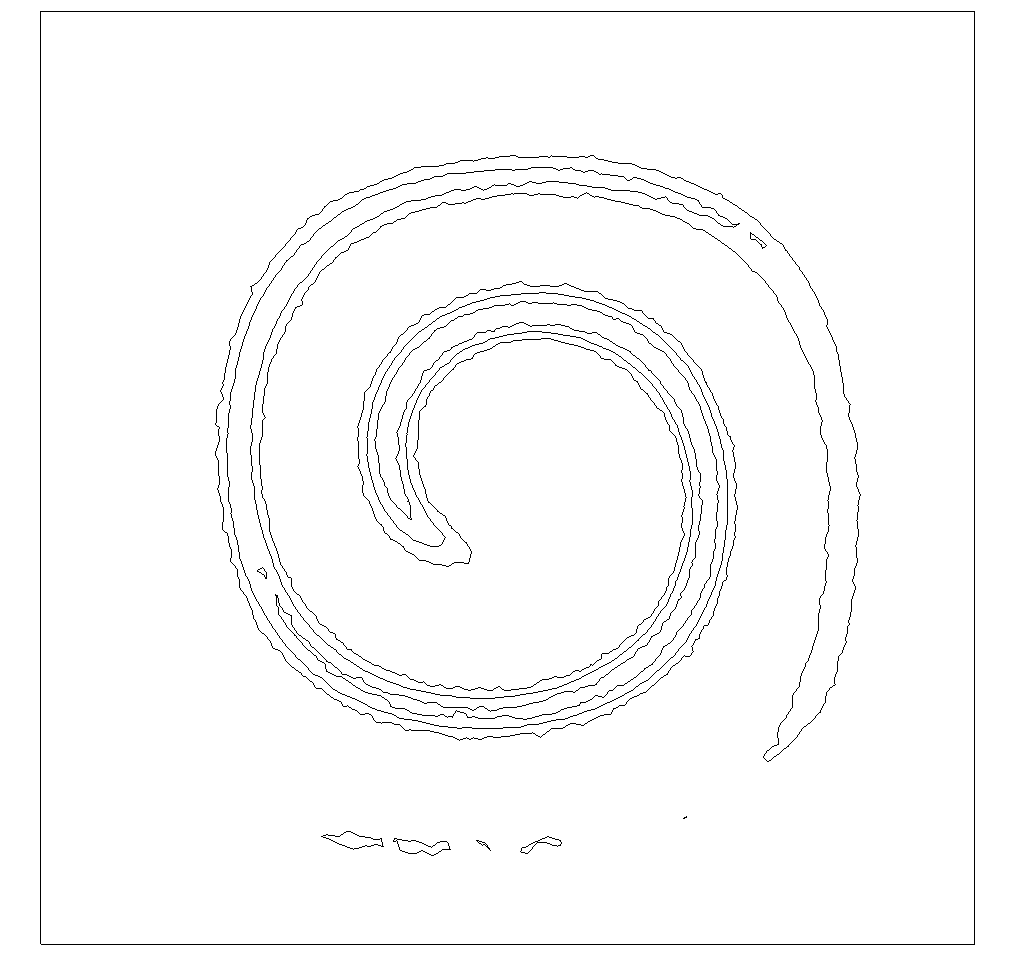}
}
\subfloat[][$\kappa_d=100$]
{
\includegraphics[width=5cm]{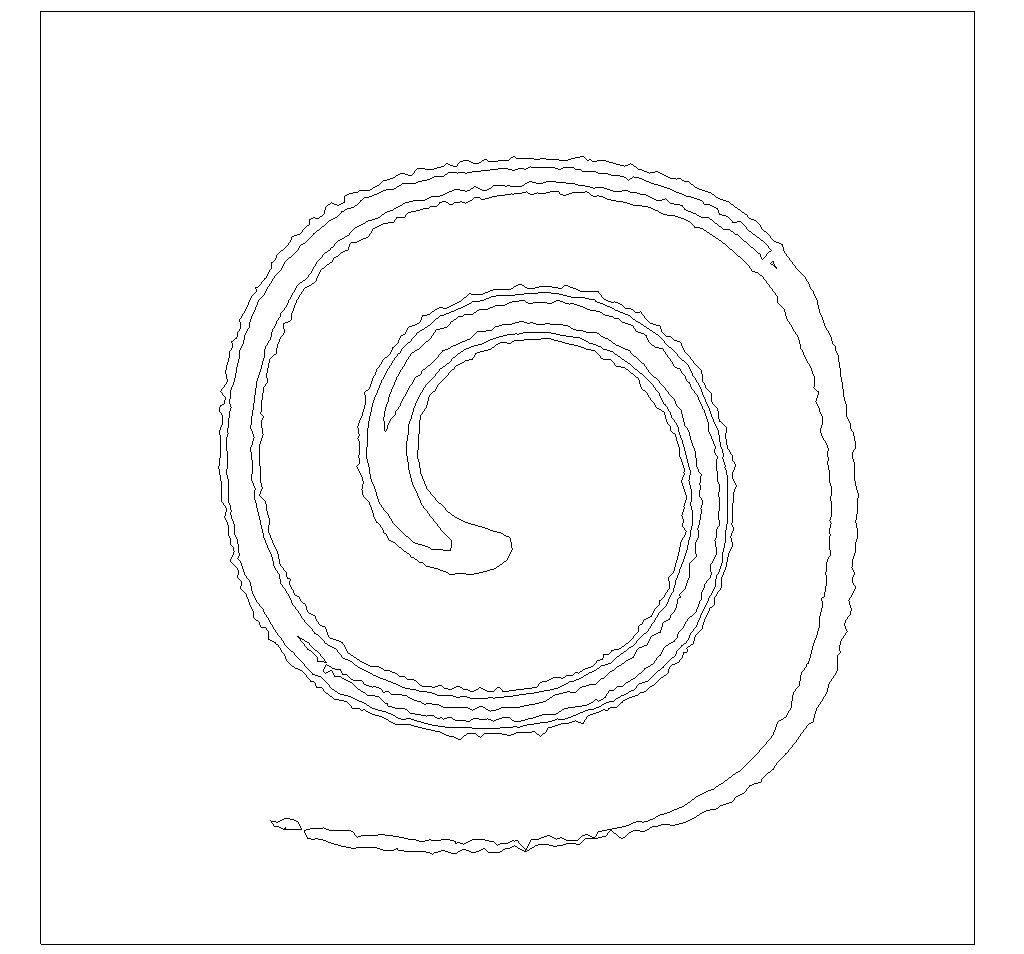}
}
\caption{Contour plots of the regularized Heaviside for {\it projected inverse scaling} on a triangular mesh with $4541$ elements.}
\label{fig:mode4_t}
\end{center}
\end{figure}

\begin{figure}[!h]
\begin{center}
\input{mode4/m4_k1_80p1_corr.tex}
\caption{Correction required for volume conservation for the $80\times80$ linear mesh with $\kappa_d=1$.}
\label{fig:mode4_corr}
\end{center}
\end{figure}

Figure \ref{fig:mode4_c1} shows the contours for the same smoothing range  but now on the $C_1-Q_2$ mesh. The results are very similar.
For the triangular mesh the results, as shown in figure \ref{fig:mode4_t} ,  are less smooth due to unstructured nature of the mesh.
Therefore, regularization has a more pronounced effect for this mesh. The higher $\kappa_d$ results in better and smoother answers.
However, as shown in section \ref{sec:dist} there is also an upper limit to the desired smoothing.

The better behaviour of the {\it projected inverse scaling} compared with the {\it projected  scaling}  
is also apparent when considering the time trace of the correction required for volume conservation, see figure \ref{fig:mode4_corr}. 
This correction is smaller, but also varies significantly less violent.

\subsection{Convergence study} 

Using the results from the previous sections we focus our attention on the approaches that showed the most potential. This 
is {\it projected inverse scaling} which is suitable for all mesh types. 
For the quadrilateral  meshes, the smoothing is absent, that is  $\kappa_d=0$, while for the triangular 
mesh we use a relatively high smoothing of $\kappa_d=10$. The previous paragraph  indicates  these smoothing values lead to satisfying results.  
 
 In this section we investigate the convergence of this approach. We  consider two different 
 measures for the accuracy of
the solution at the final time $t=8$. First, we look at the area mismatch between actual and anticipated area indicated by the level-set. 
The anticipated value is the initial area, i.e. at $t=0$.
This measure is computed as
\begin{align}
L_1(\hat{H}) = \int_{\Omega}  | H(\hat{\phi}_{_T}) - H(\hat{\phi}_{_0})| d\Omega, 
\end{align}
and is physically the most relevant  in future applications.

Additionally, we look how far the local distance field meandered from the anticipated value. Again the initial  solution is used as reference. 
For this we use the following norm
\begin{align}
L_{\infty}(\phi) =& \max_{\bx \in \Omega }\left   | \phi_T - \phi_0  \right |.
\end{align}

\begin{table}
   \centering
   \small
   \begin{tabular}{|c|cc|cc|cc|cc|cc|cc|}
      \hline
            & \multicolumn{4}{|c|}{Linear triangle}  & \multicolumn{4}{|c|}{Linear quad} & \multicolumn{4}{|c|}{Quadratic quad} \\
      \hline
       Elems     &  $L_1(\hat{H})$  & Rate     & $L_{\infty}(\phi)$ & Rate
                      &  $L_1(\hat{H})$  & Rate     & $L_{\infty}(\phi)$ & Rate
                      &  $L_1(\hat{H})$  & Rate     & $L_{\infty}(\phi)$ & Rate\\
      \hline
10		 & 1.26e-1 & 		 & 3.44e-1	 & 	    & 1.34e-1 &      & 3.88e-1 &      & 1.26e-1 &      & 3.06e-1 &      \\
20		 & 1.06e-1 & 	0.26 & 2.54e-1	 & 0.44 & 1.03e-1 & 0.38 & 2.37e-1 & 0.71 & 6.19e-2 & 1.03 & 1.49e-1 & 1.04 \\
40		 & 6.80e-2 & 	0.64	 & 1.42e-1	 & 0.84	& 5.37e-2 & 0.94 & 1.20e-1 & 0.98 & 1.89e-2 & 1.72 & 5.44e-2 & 1.45 \\
80		 & 3.20e-2 & 	1.09 & 7.95e-2	 & 0.84 & 2.52e-2 & 1.09 & 6.68e-2 & 0.85 & 7.46e-3 & 1.34 & 1.78e-2 & 1.61 \\  
       \hline
       \end{tabular}  
        \caption{Convergence of the level-set method.}
    \label{tab:convergence}        
\end{table}

\begin{figure}[h] 
\centering
\subfloat[][Mesh convergence in the $L_1(H)$-norm]
{
\scriptsize
\setlength{\unitlength}{0.240900pt}
\ifx\plotpoint\undefined\newsavebox{\plotpoint}\fi
\sbox{\plotpoint}{\rule[-0.200pt]{0.400pt}{0.400pt}}%
\begin{picture}(826,590)(0,0)
\font\gnuplot=cmr10 at 7pt
\gnuplot
\sbox{\plotpoint}{\rule[-0.200pt]{0.400pt}{0.400pt}}%
\put(134.0,92.0){\rule[-0.200pt]{2.409pt}{0.400pt}}
\put(773.0,92.0){\rule[-0.200pt]{2.409pt}{0.400pt}}
\put(134.0,119.0){\rule[-0.200pt]{2.409pt}{0.400pt}}
\put(773.0,119.0){\rule[-0.200pt]{2.409pt}{0.400pt}}
\put(134.0,141.0){\rule[-0.200pt]{2.409pt}{0.400pt}}
\put(773.0,141.0){\rule[-0.200pt]{2.409pt}{0.400pt}}
\put(134.0,159.0){\rule[-0.200pt]{2.409pt}{0.400pt}}
\put(773.0,159.0){\rule[-0.200pt]{2.409pt}{0.400pt}}
\put(134.0,175.0){\rule[-0.200pt]{2.409pt}{0.400pt}}
\put(773.0,175.0){\rule[-0.200pt]{2.409pt}{0.400pt}}
\put(134.0,189.0){\rule[-0.200pt]{2.409pt}{0.400pt}}
\put(773.0,189.0){\rule[-0.200pt]{2.409pt}{0.400pt}}
\put(134.0,202.0){\rule[-0.200pt]{4.818pt}{0.400pt}}
\put(120,202){\makebox(0,0)[r]{ 0.01}}
\put(763.0,202.0){\rule[-0.200pt]{4.818pt}{0.400pt}}
\put(134.0,285.0){\rule[-0.200pt]{2.409pt}{0.400pt}}
\put(773.0,285.0){\rule[-0.200pt]{2.409pt}{0.400pt}}
\put(134.0,334.0){\rule[-0.200pt]{2.409pt}{0.400pt}}
\put(773.0,334.0){\rule[-0.200pt]{2.409pt}{0.400pt}}
\put(134.0,368.0){\rule[-0.200pt]{2.409pt}{0.400pt}}
\put(773.0,368.0){\rule[-0.200pt]{2.409pt}{0.400pt}}
\put(134.0,395.0){\rule[-0.200pt]{2.409pt}{0.400pt}}
\put(773.0,395.0){\rule[-0.200pt]{2.409pt}{0.400pt}}
\put(134.0,417.0){\rule[-0.200pt]{2.409pt}{0.400pt}}
\put(773.0,417.0){\rule[-0.200pt]{2.409pt}{0.400pt}}
\put(134.0,435.0){\rule[-0.200pt]{2.409pt}{0.400pt}}
\put(773.0,435.0){\rule[-0.200pt]{2.409pt}{0.400pt}}
\put(134.0,451.0){\rule[-0.200pt]{2.409pt}{0.400pt}}
\put(773.0,451.0){\rule[-0.200pt]{2.409pt}{0.400pt}}
\put(134.0,465.0){\rule[-0.200pt]{2.409pt}{0.400pt}}
\put(773.0,465.0){\rule[-0.200pt]{2.409pt}{0.400pt}}
\put(134.0,478.0){\rule[-0.200pt]{4.818pt}{0.400pt}}
\put(120,478){\makebox(0,0)[r]{ 0.1}}
\put(763.0,478.0){\rule[-0.200pt]{4.818pt}{0.400pt}}
\put(134.0,561.0){\rule[-0.200pt]{2.409pt}{0.400pt}}
\put(773.0,561.0){\rule[-0.200pt]{2.409pt}{0.400pt}}
\put(164.0,92.0){\rule[-0.200pt]{0.400pt}{4.818pt}}
\put(164,63){\makebox(0,0){ 10}}
\put(164.0,541.0){\rule[-0.200pt]{0.400pt}{4.818pt}}
\put(359.0,92.0){\rule[-0.200pt]{0.400pt}{4.818pt}}
\put(359,63){\makebox(0,0){ 20}}
\put(359.0,541.0){\rule[-0.200pt]{0.400pt}{4.818pt}}
\put(554.0,92.0){\rule[-0.200pt]{0.400pt}{4.818pt}}
\put(554,63){\makebox(0,0){ 40}}
\put(554.0,541.0){\rule[-0.200pt]{0.400pt}{4.818pt}}
\put(750.0,92.0){\rule[-0.200pt]{0.400pt}{4.818pt}}
\put(750,63){\makebox(0,0){ 80}}
\put(750.0,541.0){\rule[-0.200pt]{0.400pt}{4.818pt}}
\put(134.0,92.0){\rule[-0.200pt]{0.400pt}{112.982pt}}
\put(134.0,92.0){\rule[-0.200pt]{156.344pt}{0.400pt}}
\put(783.0,92.0){\rule[-0.200pt]{0.400pt}{112.982pt}}
\put(134.0,561.0){\rule[-0.200pt]{156.344pt}{0.400pt}}
\put(21,326){\makebox(0,0){\rotatebox{90}{$L_1(H)$}}}
\put(458,20){\makebox(0,0){$N_{elem}$}}
\put(665,527){\makebox(0,0)[r]{$C_0-P_1$}}
\put(679.0,527.0){\rule[-0.200pt]{18.308pt}{0.400pt}}
\put(164,506){\usebox{\plotpoint}}
\multiput(164.00,504.92)(4.715,-0.496){39}{\rule{3.814pt}{0.119pt}}
\multiput(164.00,505.17)(187.083,-21.000){2}{\rule{1.907pt}{0.400pt}}
\multiput(359.00,483.92)(1.848,-0.498){103}{\rule{1.572pt}{0.120pt}}
\multiput(359.00,484.17)(191.738,-53.000){2}{\rule{0.786pt}{0.400pt}}
\multiput(554.00,430.92)(1.079,-0.499){179}{\rule{0.962pt}{0.120pt}}
\multiput(554.00,431.17)(194.004,-91.000){2}{\rule{0.481pt}{0.400pt}}
\put(164,506){\makebox(0,0){$+$}}
\put(359,485){\makebox(0,0){$+$}}
\put(554,432){\makebox(0,0){$+$}}
\put(750,341){\makebox(0,0){$+$}}
\put(717,527){\makebox(0,0){$+$}}
\put(665,498){\makebox(0,0)[r]{$C_0-Q_1$}}
\multiput(679,498)(20.756,0.000){4}{\usebox{\plotpoint}}
\put(755,498){\usebox{\plotpoint}}
\put(164,513){\usebox{\plotpoint}}
\multiput(164,513)(20.498,-3.259){10}{\usebox{\plotpoint}}
\multiput(359,482)(19.237,-7.793){10}{\usebox{\plotpoint}}
\multiput(554,403)(18.862,-8.661){11}{\usebox{\plotpoint}}
\put(750,313){\usebox{\plotpoint}}
\put(164,513){\makebox(0,0){$\times$}}
\put(359,482){\makebox(0,0){$\times$}}
\put(554,403){\makebox(0,0){$\times$}}
\put(750,313){\makebox(0,0){$\times$}}
\put(717,498){\makebox(0,0){$\times$}}
\sbox{\plotpoint}{\rule[-0.400pt]{0.800pt}{0.800pt}}%
\sbox{\plotpoint}{\rule[-0.200pt]{0.400pt}{0.400pt}}%
\put(665,469){\makebox(0,0)[r]{$C_1-Q_2$}}
\sbox{\plotpoint}{\rule[-0.400pt]{0.800pt}{0.800pt}}%
\put(679.0,469.0){\rule[-0.400pt]{18.308pt}{0.800pt}}
\put(164,506){\usebox{\plotpoint}}
\multiput(164.00,504.09)(1.138,-0.501){165}{\rule{2.014pt}{0.121pt}}
\multiput(164.00,504.34)(190.820,-86.000){2}{\rule{1.007pt}{0.800pt}}
\multiput(359.00,418.09)(0.687,-0.501){277}{\rule{1.299pt}{0.121pt}}
\multiput(359.00,418.34)(192.305,-142.000){2}{\rule{0.649pt}{0.800pt}}
\multiput(554.00,276.09)(0.885,-0.501){215}{\rule{1.613pt}{0.121pt}}
\multiput(554.00,276.34)(192.653,-111.000){2}{\rule{0.806pt}{0.800pt}}
\put(164,506){\makebox(0,0){$\ast$}}
\put(359,420){\makebox(0,0){$\ast$}}
\put(554,278){\makebox(0,0){$\ast$}}
\put(750,167){\makebox(0,0){$\ast$}}
\put(717,469){\makebox(0,0){$\ast$}}
\sbox{\plotpoint}{\rule[-0.200pt]{0.400pt}{0.400pt}}%
\put(134.0,92.0){\rule[-0.200pt]{0.400pt}{112.982pt}}
\put(134.0,92.0){\rule[-0.200pt]{156.344pt}{0.400pt}}
\put(783.0,92.0){\rule[-0.200pt]{0.400pt}{112.982pt}}
\put(134.0,561.0){\rule[-0.200pt]{156.344pt}{0.400pt}}
\end{picture}  
}
\subfloat[][Mesh convergence in the $L_{\infty}(\phi)$-norm]
{
\scriptsize
\setlength{\unitlength}{0.240900pt}
\ifx\plotpoint\undefined\newsavebox{\plotpoint}\fi
\begin{picture}(826,590)(0,0)
\font\gnuplot=cmr10 at 7pt
\gnuplot
\sbox{\plotpoint}{\rule[-0.200pt]{0.400pt}{0.400pt}}%
\put(134.0,92.0){\rule[-0.200pt]{4.818pt}{0.400pt}}
\put(120,92){\makebox(0,0)[r]{ 0.01}}
\put(763.0,92.0){\rule[-0.200pt]{4.818pt}{0.400pt}}
\put(134.0,180.0){\rule[-0.200pt]{2.409pt}{0.400pt}}
\put(773.0,180.0){\rule[-0.200pt]{2.409pt}{0.400pt}}
\put(134.0,232.0){\rule[-0.200pt]{2.409pt}{0.400pt}}
\put(773.0,232.0){\rule[-0.200pt]{2.409pt}{0.400pt}}
\put(134.0,268.0){\rule[-0.200pt]{2.409pt}{0.400pt}}
\put(773.0,268.0){\rule[-0.200pt]{2.409pt}{0.400pt}}
\put(134.0,297.0){\rule[-0.200pt]{2.409pt}{0.400pt}}
\put(773.0,297.0){\rule[-0.200pt]{2.409pt}{0.400pt}}
\put(134.0,320.0){\rule[-0.200pt]{2.409pt}{0.400pt}}
\put(773.0,320.0){\rule[-0.200pt]{2.409pt}{0.400pt}}
\put(134.0,339.0){\rule[-0.200pt]{2.409pt}{0.400pt}}
\put(773.0,339.0){\rule[-0.200pt]{2.409pt}{0.400pt}}
\put(134.0,356.0){\rule[-0.200pt]{2.409pt}{0.400pt}}
\put(773.0,356.0){\rule[-0.200pt]{2.409pt}{0.400pt}}
\put(134.0,371.0){\rule[-0.200pt]{2.409pt}{0.400pt}}
\put(773.0,371.0){\rule[-0.200pt]{2.409pt}{0.400pt}}
\put(134.0,385.0){\rule[-0.200pt]{4.818pt}{0.400pt}}
\put(120,385){\makebox(0,0)[r]{ 0.1}}
\put(763.0,385.0){\rule[-0.200pt]{4.818pt}{0.400pt}}
\put(134.0,473.0){\rule[-0.200pt]{2.409pt}{0.400pt}}
\put(773.0,473.0){\rule[-0.200pt]{2.409pt}{0.400pt}}
\put(134.0,524.0){\rule[-0.200pt]{2.409pt}{0.400pt}}
\put(773.0,524.0){\rule[-0.200pt]{2.409pt}{0.400pt}}
\put(134.0,561.0){\rule[-0.200pt]{2.409pt}{0.400pt}}
\put(773.0,561.0){\rule[-0.200pt]{2.409pt}{0.400pt}}
\put(164.0,92.0){\rule[-0.200pt]{0.400pt}{4.818pt}}
\put(164,63){\makebox(0,0){ 10}}
\put(164.0,541.0){\rule[-0.200pt]{0.400pt}{4.818pt}}
\put(359.0,92.0){\rule[-0.200pt]{0.400pt}{4.818pt}}
\put(359,63){\makebox(0,0){ 20}}
\put(359.0,541.0){\rule[-0.200pt]{0.400pt}{4.818pt}}
\put(554.0,92.0){\rule[-0.200pt]{0.400pt}{4.818pt}}
\put(554,63){\makebox(0,0){ 40}}
\put(554.0,541.0){\rule[-0.200pt]{0.400pt}{4.818pt}}
\put(750.0,92.0){\rule[-0.200pt]{0.400pt}{4.818pt}}
\put(750,63){\makebox(0,0){ 80}}
\put(750.0,541.0){\rule[-0.200pt]{0.400pt}{4.818pt}}
\put(134.0,92.0){\rule[-0.200pt]{0.400pt}{112.982pt}}
\put(134.0,92.0){\rule[-0.200pt]{156.344pt}{0.400pt}}
\put(783.0,92.0){\rule[-0.200pt]{0.400pt}{112.982pt}}
\put(134.0,561.0){\rule[-0.200pt]{156.344pt}{0.400pt}}
\put(21,326){\makebox(0,0){\rotatebox{90}{$L_{\infty}(\phi)$}}}
\put(458,20){\makebox(0,0){$N_{elem}$}}
\put(665,527){\makebox(0,0)[r]{$C_0-P_1$}}
\put(679.0,527.0){\rule[-0.200pt]{18.308pt}{0.400pt}}
\put(164,542){\usebox{\plotpoint}}
\multiput(164.00,540.92)(2.517,-0.498){75}{\rule{2.100pt}{0.120pt}}
\multiput(164.00,541.17)(190.641,-39.000){2}{\rule{1.050pt}{0.400pt}}
\multiput(359.00,501.92)(1.339,-0.499){143}{\rule{1.168pt}{0.120pt}}
\multiput(359.00,502.17)(192.575,-73.000){2}{\rule{0.584pt}{0.400pt}}
\multiput(554.00,428.92)(1.328,-0.499){145}{\rule{1.159pt}{0.120pt}}
\multiput(554.00,429.17)(193.593,-74.000){2}{\rule{0.580pt}{0.400pt}}
\put(164,542){\makebox(0,0){$+$}}
\put(359,503){\makebox(0,0){$+$}}
\put(554,430){\makebox(0,0){$+$}}
\put(750,356){\makebox(0,0){$+$}}
\put(717,527){\makebox(0,0){$+$}}
\put(665,498){\makebox(0,0)[r]{$C_0-Q_1$}}
\multiput(679,498)(20.756,0.000){4}{\usebox{\plotpoint}}
\put(755,498){\usebox{\plotpoint}}
\put(164,557){\usebox{\plotpoint}}
\multiput(164,557)(19.750,-6.381){10}{\usebox{\plotpoint}}
\multiput(359,494)(18.991,-8.375){11}{\usebox{\plotpoint}}
\multiput(554,408)(19.385,-7.418){10}{\usebox{\plotpoint}}
\put(750,333){\usebox{\plotpoint}}
\put(164,557){\makebox(0,0){$\times$}}
\put(359,494){\makebox(0,0){$\times$}}
\put(554,408){\makebox(0,0){$\times$}}
\put(750,333){\makebox(0,0){$\times$}}
\put(717,498){\makebox(0,0){$\times$}}
\sbox{\plotpoint}{\rule[-0.400pt]{0.800pt}{0.800pt}}%
\sbox{\plotpoint}{\rule[-0.200pt]{0.400pt}{0.400pt}}%
\put(665,469){\makebox(0,0)[r]{$C_1-Q_2$}}
\sbox{\plotpoint}{\rule[-0.400pt]{0.800pt}{0.800pt}}%
\put(679.0,469.0){\rule[-0.400pt]{18.308pt}{0.800pt}}
\put(164,527){\usebox{\plotpoint}}
\multiput(164.00,525.09)(1.063,-0.501){177}{\rule{1.896pt}{0.121pt}}
\multiput(164.00,525.34)(191.065,-92.000){2}{\rule{0.948pt}{0.800pt}}
\multiput(359.00,433.09)(0.763,-0.501){249}{\rule{1.419pt}{0.121pt}}
\multiput(359.00,433.34)(192.055,-128.000){2}{\rule{0.709pt}{0.800pt}}
\multiput(554.00,305.09)(0.691,-0.501){277}{\rule{1.304pt}{0.121pt}}
\multiput(554.00,305.34)(193.293,-142.000){2}{\rule{0.652pt}{0.800pt}}
\put(164,527){\makebox(0,0){$\ast$}}
\put(359,435){\makebox(0,0){$\ast$}}
\put(554,307){\makebox(0,0){$\ast$}}
\put(750,165){\makebox(0,0){$\ast$}}
\put(717,469){\makebox(0,0){$\ast$}}
\sbox{\plotpoint}{\rule[-0.200pt]{0.400pt}{0.400pt}}%
\put(134.0,92.0){\rule[-0.200pt]{0.400pt}{112.982pt}}
\put(134.0,92.0){\rule[-0.200pt]{156.344pt}{0.400pt}}
\put(783.0,92.0){\rule[-0.200pt]{0.400pt}{112.982pt}}
\put(134.0,561.0){\rule[-0.200pt]{156.344pt}{0.400pt}}
\end{picture}
}
\caption{Mesh convergences plots of the {\it projected inverse scaling} approach. }
\label{fig:conv}
\end{figure}

The convergence data for a sequence of meshes is given in table \ref{tab:convergence}  and is plotted in  
figure  \ref{fig:conv}. 
For the triangular mesh the element count used in table \ref{tab:convergence} and figure  \ref{fig:conv}  is the target line elements along each of the edges of the square.
The resulting meshes have $242$,  $1054$,  $4262$  and $16794$ triangular elements and
$141$, $568$,  $2212$ and $8558$ nodes, respectively. 

The convergence results are less then expected. 
For the linear meshes the convergence  approaches first order for both norms, while for the 
quadratic NURBS mesh both norms converge with a slightly higher rate of around $1.5$. 
This is less than is reported in \cite{ABKF11}. 
An explanation could be the difference in volume conservation methodology:
local  volume conservation is used \cite{ABKF11}, while here
global volume conservation is employed.

\subsection{Three-dimensional vortex in a cube}
\label{sec:3D}

In this section we investigate the applicability of the proposed method in a 3D setting. 
For this we adopt a 3D version of the vortex in box problem \cite{HighResLeveque}.
This results in the following imposed velocity field:
\begin{align}
u =& 2 \cos \left (\frac{\pi t}{3} \right ) \sin(\pi x)^2 \sin (2\pi y) \sin (2\pi z),   \nonumber \\
 v =&  -\cos \left (\frac{\pi t}{3} \right ) \sin (2\pi x) \sin(\pi y)^2 \sin (2\pi z),   \nonumber\\
  v =&  -\cos \left (\frac{\pi t}{3} \right ) \sin (2\pi x) \sin(2 \pi y) \sin (\pi z)^2.
\end{align}
The level-set experiences a full cycle of deformation in $T=3.0$, similar as in  \cite{FLD1475}. 
The 3D deformation is the superimposed two deformations used in the 2D cases, as given in eq ( \ref{eq:2Ddef}). Due to this extreme deformation the volume is spread extremely thin. This is a 
hard problem to resolve, therefore, we use a mesh with $128\times128\times128$ linear NURBS elements.  Similar to the convergence study we only consider {\it projected inverse scaling}.

\begin{figure}[!h]
\begin{center}
\subfloat[][$t=0.0$]
{
\includegraphics[width=3.8cm]{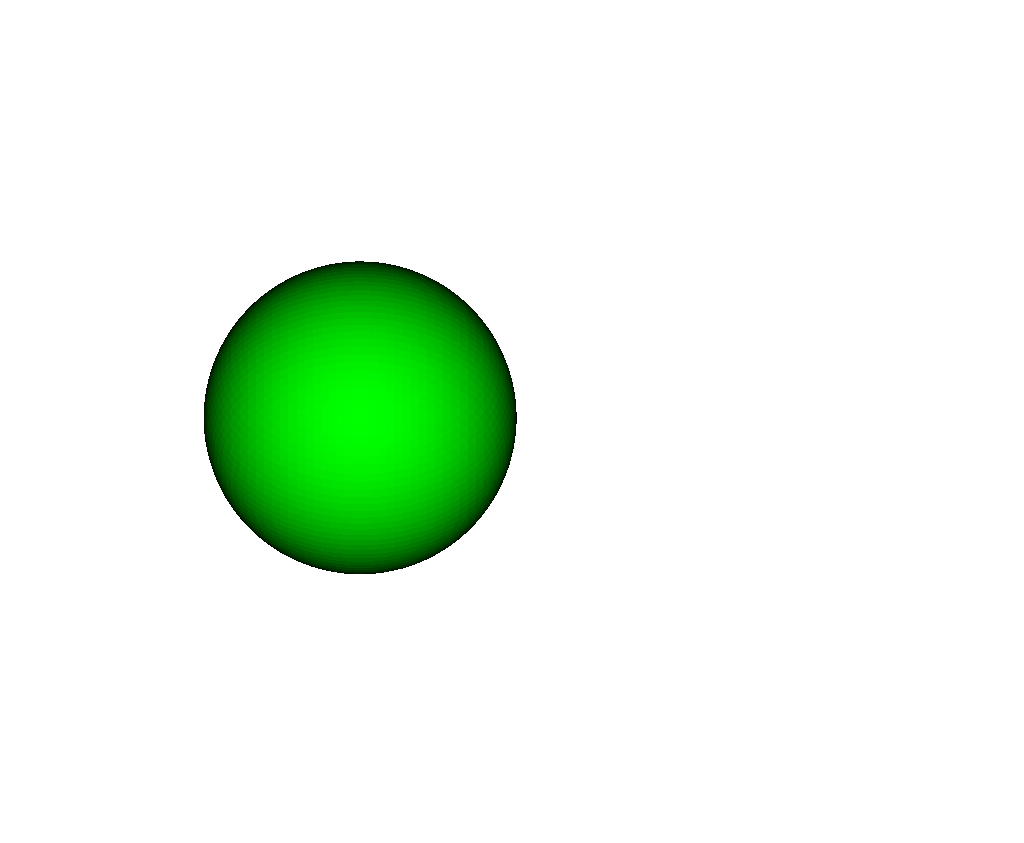}
}
\subfloat[][$t=0.2$]
{
\includegraphics[width=3.8cm]{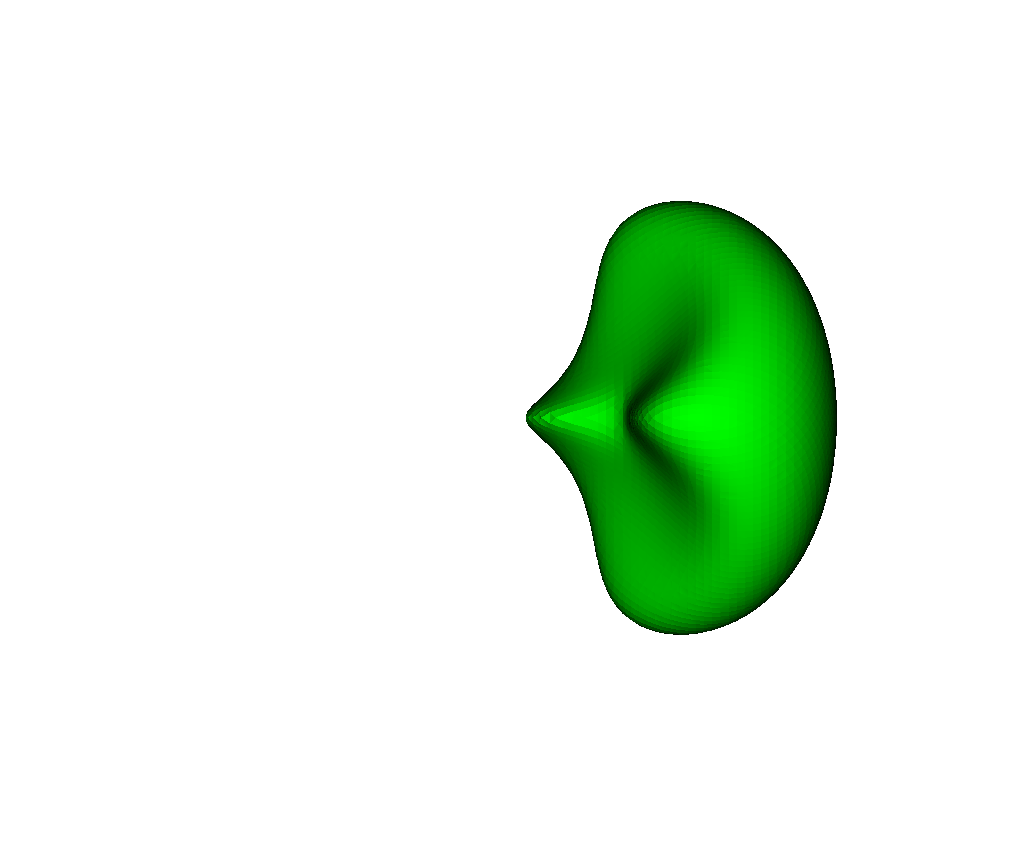}
}
\subfloat[][$t=0.4$]
{
\includegraphics[width=3.8cm]{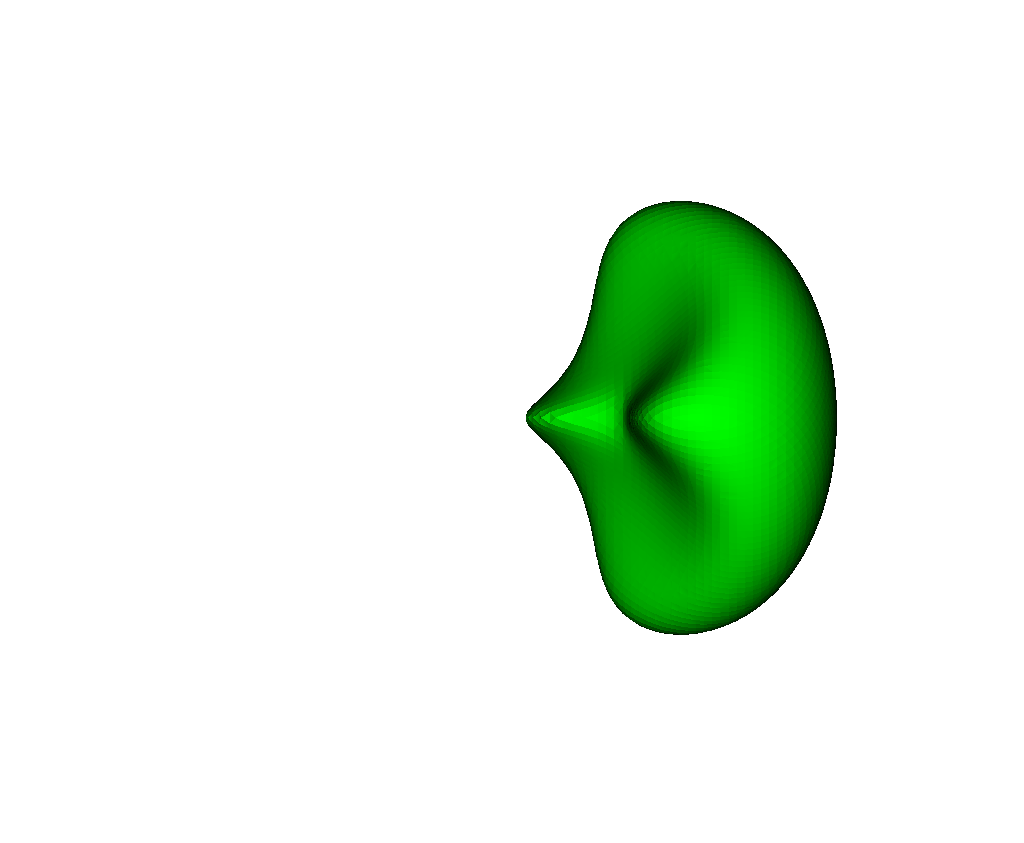}
}
\subfloat[][$t=0.6$]
{
\includegraphics[width=3.8cm]{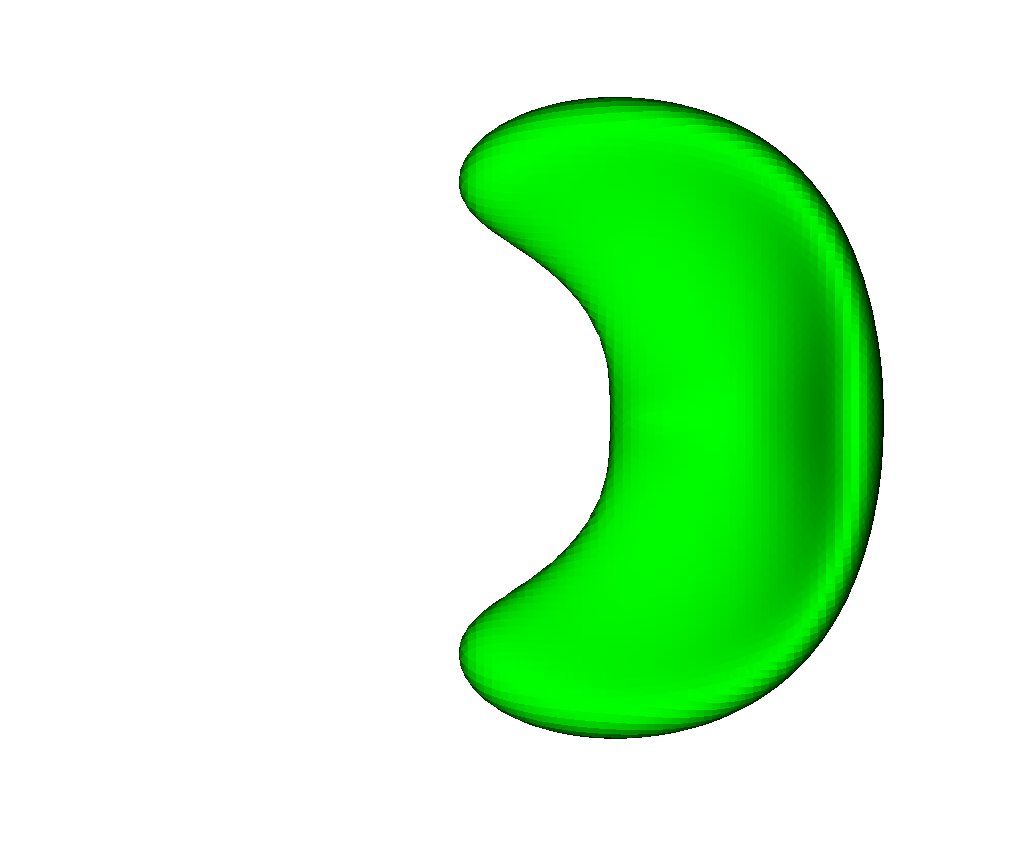}
}
\\
\subfloat[][$t=1.4$]
{
\includegraphics[width=3.8cm]{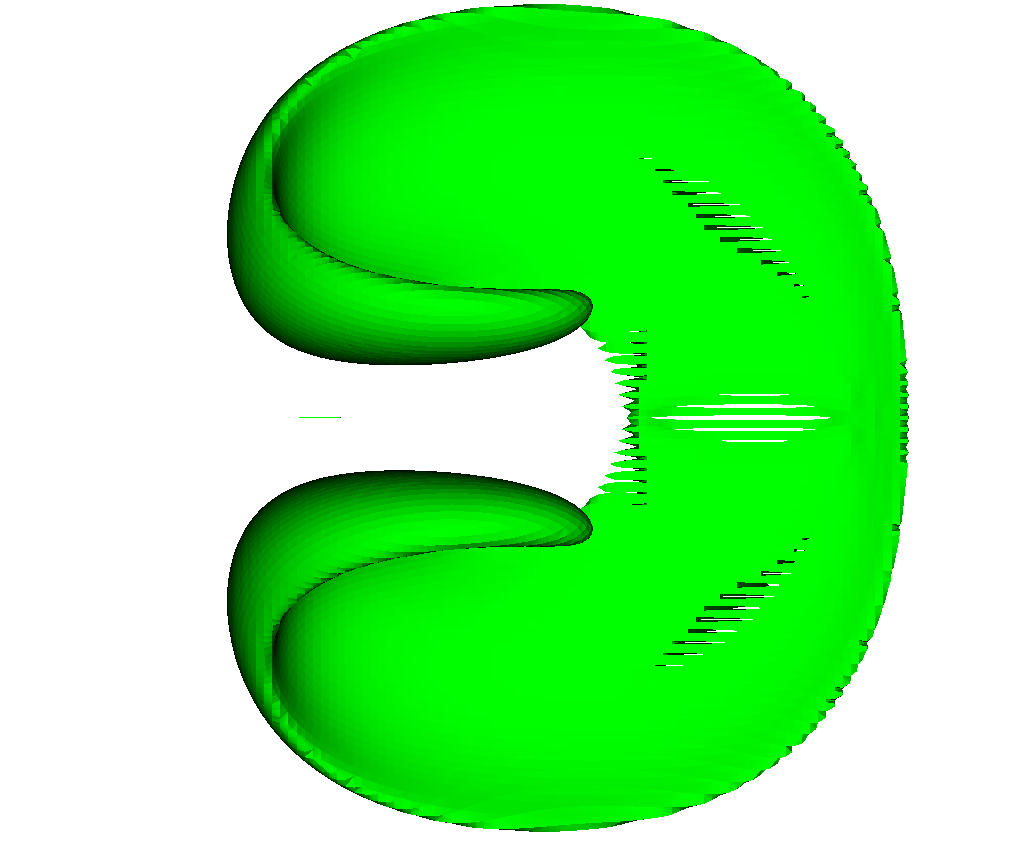}
}
\subfloat[][$t=1.8$]
{
\includegraphics[width=3.8cm]{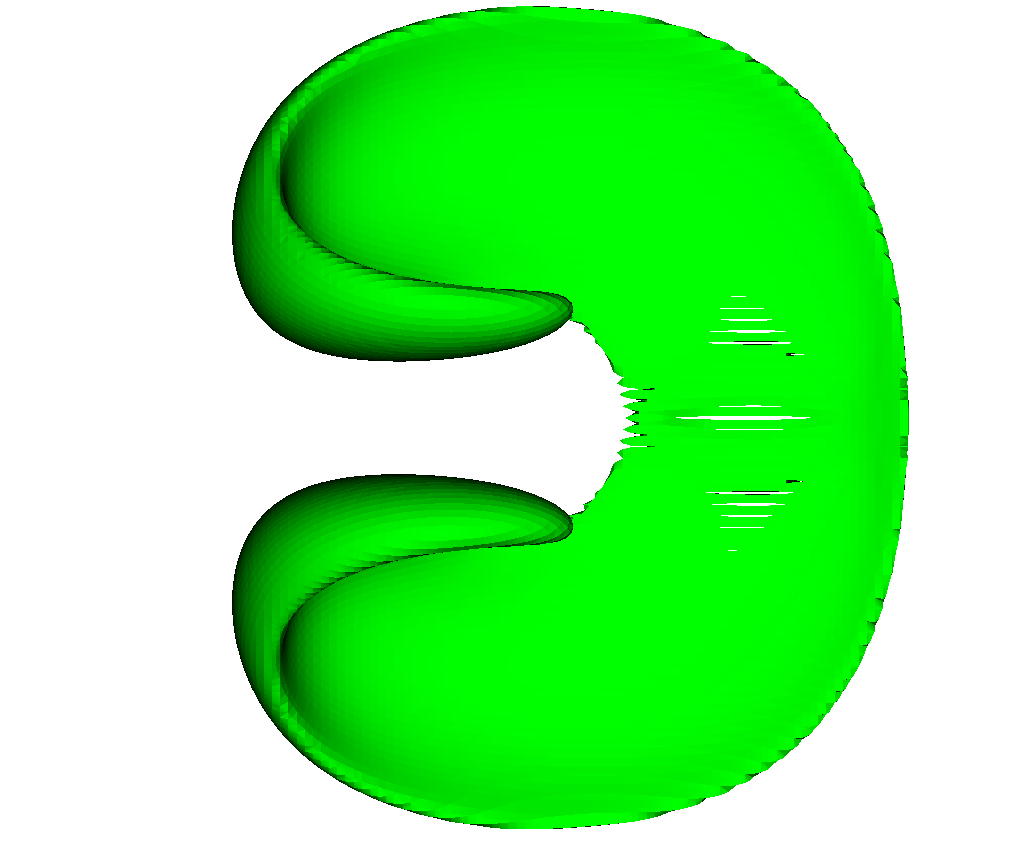}
}
\subfloat[][$t=2.6$]
{
\includegraphics[width=3.8cm]{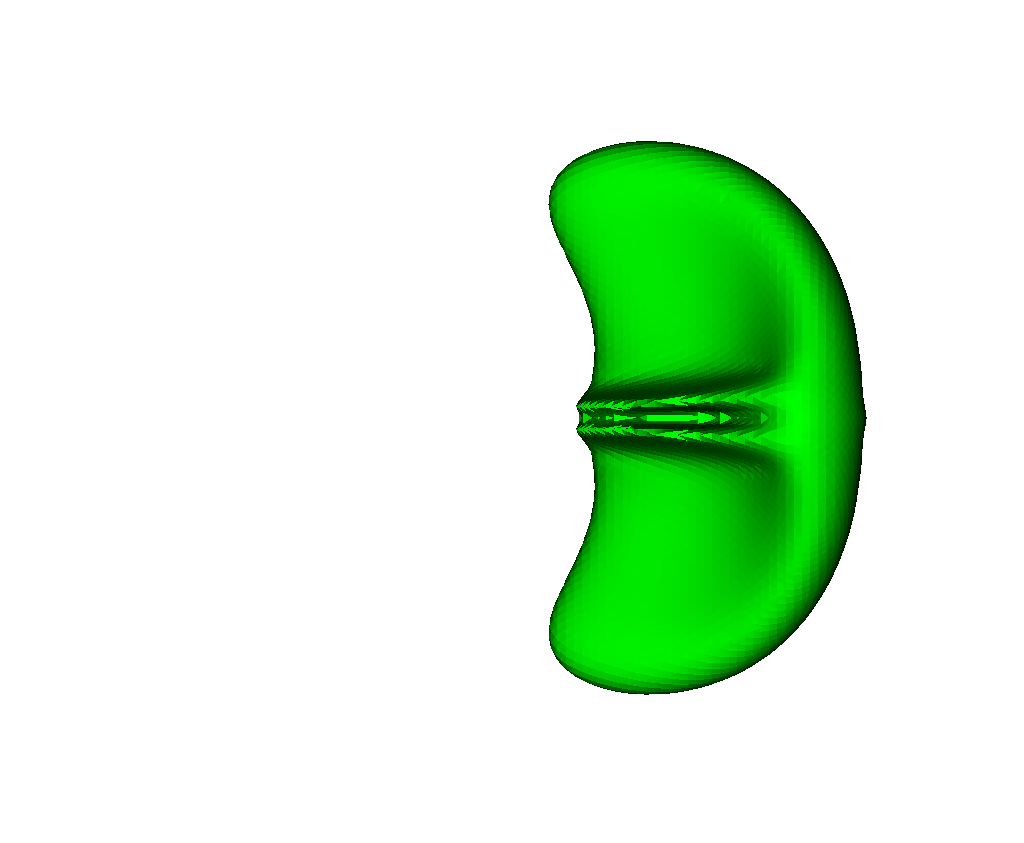}
}
\subfloat[][$t=3.0$]
{
\includegraphics[width=3.8cm]{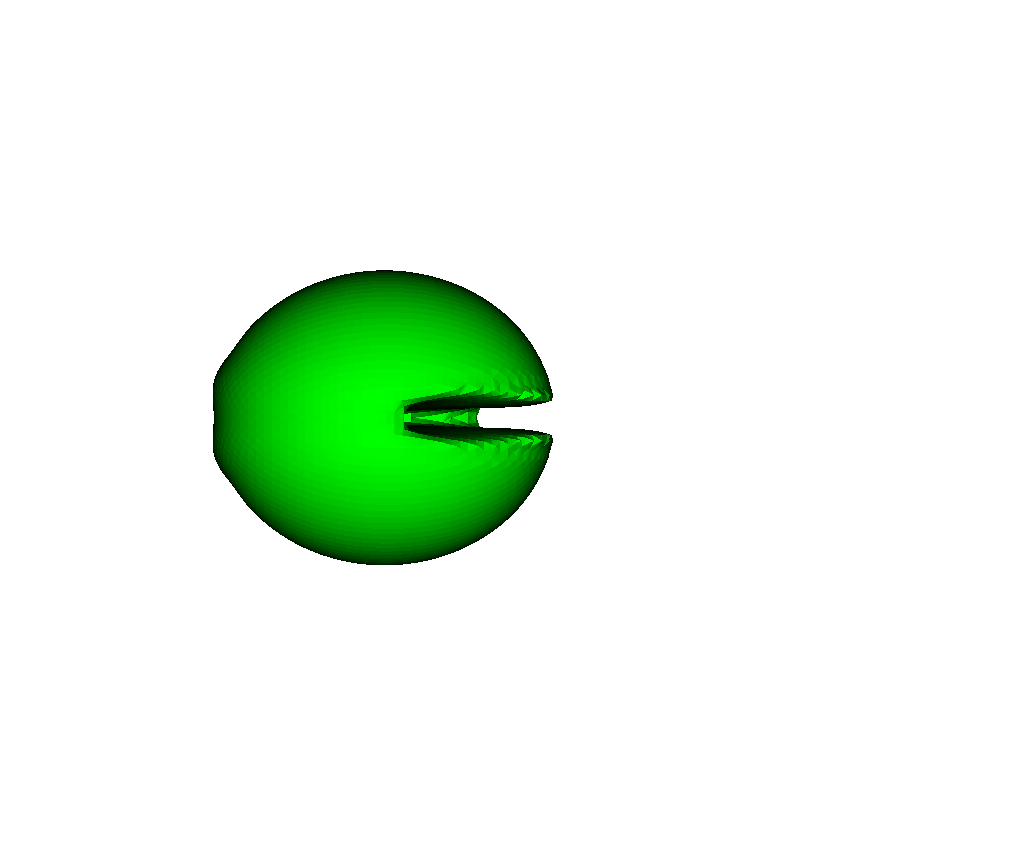}
}
\caption{Three-dimensional vortex snapshots on a mesh with $128^3$ linear NURBS elements.}
\label{fig:3Dsnap}
\end{center}
\end{figure}

In figure \ref{fig:3Dsnap} eight snapshots are of the $\phi=0$ level-set are shown.  This level-set 
represents the middle of the interface region. 
The snapshots at $t=1.4$ and $t=1.8$ the volume gets distorted to a high degree, resulting in very thin filaments.

\begin{figure}[!ht]
\begin{center}
\subfloat[][Volume error ]{
\setlength{\unitlength}{0.240900pt}
\ifx\plotpoint\undefined\newsavebox{\plotpoint}\fi
\sbox{\plotpoint}{\rule[-0.200pt]{0.400pt}{0.400pt}}%
\begin{picture}(944,590)(0,0)
\sbox{\plotpoint}{\rule[-0.200pt]{0.400pt}{0.400pt}}%
\put(211.0,131.0){\rule[-0.200pt]{4.818pt}{0.400pt}}
\put(191,131){\makebox(0,0)[r]{-0.025}}
\put(863.0,131.0){\rule[-0.200pt]{4.818pt}{0.400pt}}
\put(211.0,201.0){\rule[-0.200pt]{4.818pt}{0.400pt}}
\put(191,201){\makebox(0,0)[r]{-0.02}}
\put(863.0,201.0){\rule[-0.200pt]{4.818pt}{0.400pt}}
\put(211.0,270.0){\rule[-0.200pt]{4.818pt}{0.400pt}}
\put(191,270){\makebox(0,0)[r]{-0.015}}
\put(863.0,270.0){\rule[-0.200pt]{4.818pt}{0.400pt}}
\put(211.0,340.0){\rule[-0.200pt]{4.818pt}{0.400pt}}
\put(191,340){\makebox(0,0)[r]{-0.01}}
\put(863.0,340.0){\rule[-0.200pt]{4.818pt}{0.400pt}}
\put(211.0,410.0){\rule[-0.200pt]{4.818pt}{0.400pt}}
\put(191,410){\makebox(0,0)[r]{-0.005}}
\put(863.0,410.0){\rule[-0.200pt]{4.818pt}{0.400pt}}
\put(211.0,479.0){\rule[-0.200pt]{4.818pt}{0.400pt}}
\put(191,479){\makebox(0,0)[r]{ 0}}
\put(863.0,479.0){\rule[-0.200pt]{4.818pt}{0.400pt}}
\put(211.0,549.0){\rule[-0.200pt]{4.818pt}{0.400pt}}
\put(191,549){\makebox(0,0)[r]{ 0.005}}
\put(863.0,549.0){\rule[-0.200pt]{4.818pt}{0.400pt}}
\put(211.0,131.0){\rule[-0.200pt]{0.400pt}{4.818pt}}
\put(211,90){\makebox(0,0){ 0}}
\put(211.0,529.0){\rule[-0.200pt]{0.400pt}{4.818pt}}
\put(323.0,131.0){\rule[-0.200pt]{0.400pt}{4.818pt}}
\put(323,90){\makebox(0,0){ 0.5}}
\put(323.0,529.0){\rule[-0.200pt]{0.400pt}{4.818pt}}
\put(435.0,131.0){\rule[-0.200pt]{0.400pt}{4.818pt}}
\put(435,90){\makebox(0,0){ 1}}
\put(435.0,529.0){\rule[-0.200pt]{0.400pt}{4.818pt}}
\put(547.0,131.0){\rule[-0.200pt]{0.400pt}{4.818pt}}
\put(547,90){\makebox(0,0){ 1.5}}
\put(547.0,529.0){\rule[-0.200pt]{0.400pt}{4.818pt}}
\put(659.0,131.0){\rule[-0.200pt]{0.400pt}{4.818pt}}
\put(659,90){\makebox(0,0){ 2}}
\put(659.0,529.0){\rule[-0.200pt]{0.400pt}{4.818pt}}
\put(771.0,131.0){\rule[-0.200pt]{0.400pt}{4.818pt}}
\put(771,90){\makebox(0,0){ 2.5}}
\put(771.0,529.0){\rule[-0.200pt]{0.400pt}{4.818pt}}
\put(883.0,131.0){\rule[-0.200pt]{0.400pt}{4.818pt}}
\put(883,90){\makebox(0,0){ 3}}
\put(883.0,529.0){\rule[-0.200pt]{0.400pt}{4.818pt}}
\put(211.0,131.0){\rule[-0.200pt]{0.400pt}{100.696pt}}
\put(211.0,131.0){\rule[-0.200pt]{161.885pt}{0.400pt}}
\put(883.0,131.0){\rule[-0.200pt]{0.400pt}{100.696pt}}
\put(211.0,549.0){\rule[-0.200pt]{161.885pt}{0.400pt}}
\put(30,340){\makebox(0,0){\rotatebox{90}{$\Delta V~~(\%)$}}}
\put(547,29){\makebox(0,0){$t$}}
\put(211,479){\usebox{\plotpoint}}
\put(211,477.67){\rule{2.650pt}{0.400pt}}
\multiput(211.00,478.17)(5.500,-1.000){2}{\rule{1.325pt}{0.400pt}}
\put(267,477.67){\rule{2.650pt}{0.400pt}}
\multiput(267.00,477.17)(5.500,1.000){2}{\rule{1.325pt}{0.400pt}}
\put(222.0,478.0){\rule[-0.200pt]{10.840pt}{0.400pt}}
\put(301,477.17){\rule{2.300pt}{0.400pt}}
\multiput(301.00,478.17)(6.226,-2.000){2}{\rule{1.150pt}{0.400pt}}
\multiput(312.00,475.94)(1.505,-0.468){5}{\rule{1.200pt}{0.113pt}}
\multiput(312.00,476.17)(8.509,-4.000){2}{\rule{0.600pt}{0.400pt}}
\multiput(323.00,473.59)(1.155,0.477){7}{\rule{0.980pt}{0.115pt}}
\multiput(323.00,472.17)(8.966,5.000){2}{\rule{0.490pt}{0.400pt}}
\multiput(334.00,478.61)(2.248,0.447){3}{\rule{1.567pt}{0.108pt}}
\multiput(334.00,477.17)(7.748,3.000){2}{\rule{0.783pt}{0.400pt}}
\multiput(345.58,477.96)(0.492,-0.798){21}{\rule{0.119pt}{0.733pt}}
\multiput(344.17,479.48)(12.000,-17.478){2}{\rule{0.400pt}{0.367pt}}
\multiput(357.00,460.93)(1.155,-0.477){7}{\rule{0.980pt}{0.115pt}}
\multiput(357.00,461.17)(8.966,-5.000){2}{\rule{0.490pt}{0.400pt}}
\multiput(368.58,452.51)(0.492,-1.251){19}{\rule{0.118pt}{1.082pt}}
\multiput(367.17,454.75)(11.000,-24.755){2}{\rule{0.400pt}{0.541pt}}
\multiput(379.58,423.09)(0.492,-2.005){19}{\rule{0.118pt}{1.664pt}}
\multiput(378.17,426.55)(11.000,-39.547){2}{\rule{0.400pt}{0.832pt}}
\multiput(390.58,376.62)(0.492,-3.090){19}{\rule{0.118pt}{2.500pt}}
\multiput(389.17,381.81)(11.000,-60.811){2}{\rule{0.400pt}{1.250pt}}
\multiput(401.58,311.59)(0.492,-2.780){21}{\rule{0.119pt}{2.267pt}}
\multiput(400.17,316.30)(12.000,-60.295){2}{\rule{0.400pt}{1.133pt}}
\multiput(413.58,250.45)(0.492,-1.581){19}{\rule{0.118pt}{1.336pt}}
\multiput(412.17,253.23)(11.000,-31.226){2}{\rule{0.400pt}{0.668pt}}
\multiput(424.58,217.51)(0.492,-1.251){19}{\rule{0.118pt}{1.082pt}}
\multiput(423.17,219.75)(11.000,-24.755){2}{\rule{0.400pt}{0.541pt}}
\multiput(435.00,193.95)(2.248,-0.447){3}{\rule{1.567pt}{0.108pt}}
\multiput(435.00,194.17)(7.748,-3.000){2}{\rule{0.783pt}{0.400pt}}
\multiput(446.58,189.62)(0.492,-0.590){19}{\rule{0.118pt}{0.573pt}}
\multiput(445.17,190.81)(11.000,-11.811){2}{\rule{0.400pt}{0.286pt}}
\put(278.0,479.0){\rule[-0.200pt]{5.541pt}{0.400pt}}
\multiput(469.00,177.93)(0.798,-0.485){11}{\rule{0.729pt}{0.117pt}}
\multiput(469.00,178.17)(9.488,-7.000){2}{\rule{0.364pt}{0.400pt}}
\multiput(480.00,170.93)(0.943,-0.482){9}{\rule{0.833pt}{0.116pt}}
\multiput(480.00,171.17)(9.270,-6.000){2}{\rule{0.417pt}{0.400pt}}
\put(457.0,179.0){\rule[-0.200pt]{2.891pt}{0.400pt}}
\put(502,164.67){\rule{2.650pt}{0.400pt}}
\multiput(502.00,165.17)(5.500,-1.000){2}{\rule{1.325pt}{0.400pt}}
\multiput(513.00,163.93)(1.267,-0.477){7}{\rule{1.060pt}{0.115pt}}
\multiput(513.00,164.17)(9.800,-5.000){2}{\rule{0.530pt}{0.400pt}}
\multiput(525.00,158.93)(1.155,-0.477){7}{\rule{0.980pt}{0.115pt}}
\multiput(525.00,159.17)(8.966,-5.000){2}{\rule{0.490pt}{0.400pt}}
\multiput(536.00,153.95)(2.248,-0.447){3}{\rule{1.567pt}{0.108pt}}
\multiput(536.00,154.17)(7.748,-3.000){2}{\rule{0.783pt}{0.400pt}}
\put(547,152.17){\rule{2.300pt}{0.400pt}}
\multiput(547.00,151.17)(6.226,2.000){2}{\rule{1.150pt}{0.400pt}}
\multiput(558.00,154.59)(0.611,0.489){15}{\rule{0.589pt}{0.118pt}}
\multiput(558.00,153.17)(9.778,9.000){2}{\rule{0.294pt}{0.400pt}}
\multiput(569.00,163.58)(0.496,0.492){21}{\rule{0.500pt}{0.119pt}}
\multiput(569.00,162.17)(10.962,12.000){2}{\rule{0.250pt}{0.400pt}}
\multiput(581.58,175.00)(0.492,0.590){19}{\rule{0.118pt}{0.573pt}}
\multiput(580.17,175.00)(11.000,11.811){2}{\rule{0.400pt}{0.286pt}}
\multiput(592.00,188.58)(0.496,0.492){19}{\rule{0.500pt}{0.118pt}}
\multiput(592.00,187.17)(9.962,11.000){2}{\rule{0.250pt}{0.400pt}}
\multiput(603.58,199.00)(0.492,0.920){19}{\rule{0.118pt}{0.827pt}}
\multiput(602.17,199.00)(11.000,18.283){2}{\rule{0.400pt}{0.414pt}}
\multiput(614.58,219.00)(0.492,1.345){19}{\rule{0.118pt}{1.155pt}}
\multiput(613.17,219.00)(11.000,26.604){2}{\rule{0.400pt}{0.577pt}}
\multiput(625.58,248.00)(0.492,0.755){21}{\rule{0.119pt}{0.700pt}}
\multiput(624.17,248.00)(12.000,16.547){2}{\rule{0.400pt}{0.350pt}}
\multiput(637.58,266.00)(0.492,1.062){19}{\rule{0.118pt}{0.936pt}}
\multiput(636.17,266.00)(11.000,21.057){2}{\rule{0.400pt}{0.468pt}}
\multiput(648.58,289.00)(0.492,0.873){19}{\rule{0.118pt}{0.791pt}}
\multiput(647.17,289.00)(11.000,17.358){2}{\rule{0.400pt}{0.395pt}}
\multiput(659.58,308.00)(0.492,2.005){19}{\rule{0.118pt}{1.664pt}}
\multiput(658.17,308.00)(11.000,39.547){2}{\rule{0.400pt}{0.832pt}}
\multiput(670.58,351.00)(0.492,1.251){19}{\rule{0.118pt}{1.082pt}}
\multiput(669.17,351.00)(11.000,24.755){2}{\rule{0.400pt}{0.541pt}}
\multiput(681.58,378.00)(0.492,1.444){21}{\rule{0.119pt}{1.233pt}}
\multiput(680.17,378.00)(12.000,31.440){2}{\rule{0.400pt}{0.617pt}}
\multiput(693.58,412.00)(0.492,1.392){19}{\rule{0.118pt}{1.191pt}}
\multiput(692.17,412.00)(11.000,27.528){2}{\rule{0.400pt}{0.595pt}}
\multiput(704.58,442.00)(0.492,0.732){19}{\rule{0.118pt}{0.682pt}}
\multiput(703.17,442.00)(11.000,14.585){2}{\rule{0.400pt}{0.341pt}}
\multiput(715.00,458.61)(2.248,0.447){3}{\rule{1.567pt}{0.108pt}}
\multiput(715.00,457.17)(7.748,3.000){2}{\rule{0.783pt}{0.400pt}}
\multiput(726.00,461.58)(0.496,0.492){19}{\rule{0.500pt}{0.118pt}}
\multiput(726.00,460.17)(9.962,11.000){2}{\rule{0.250pt}{0.400pt}}
\multiput(737.00,472.59)(0.758,0.488){13}{\rule{0.700pt}{0.117pt}}
\multiput(737.00,471.17)(10.547,8.000){2}{\rule{0.350pt}{0.400pt}}
\multiput(749.00,478.95)(2.248,-0.447){3}{\rule{1.567pt}{0.108pt}}
\multiput(749.00,479.17)(7.748,-3.000){2}{\rule{0.783pt}{0.400pt}}
\put(760,475.67){\rule{2.650pt}{0.400pt}}
\multiput(760.00,476.17)(5.500,-1.000){2}{\rule{1.325pt}{0.400pt}}
\multiput(771.00,476.59)(0.943,0.482){9}{\rule{0.833pt}{0.116pt}}
\multiput(771.00,475.17)(9.270,6.000){2}{\rule{0.417pt}{0.400pt}}
\put(782,481.67){\rule{2.650pt}{0.400pt}}
\multiput(782.00,481.17)(5.500,1.000){2}{\rule{1.325pt}{0.400pt}}
\put(793,483.17){\rule{2.500pt}{0.400pt}}
\multiput(793.00,482.17)(6.811,2.000){2}{\rule{1.250pt}{0.400pt}}
\multiput(805.00,485.60)(1.505,0.468){5}{\rule{1.200pt}{0.113pt}}
\multiput(805.00,484.17)(8.509,4.000){2}{\rule{0.600pt}{0.400pt}}
\put(816,489.17){\rule{2.300pt}{0.400pt}}
\multiput(816.00,488.17)(6.226,2.000){2}{\rule{1.150pt}{0.400pt}}
\multiput(827.00,489.95)(2.248,-0.447){3}{\rule{1.567pt}{0.108pt}}
\multiput(827.00,490.17)(7.748,-3.000){2}{\rule{0.783pt}{0.400pt}}
\multiput(838.00,486.94)(1.505,-0.468){5}{\rule{1.200pt}{0.113pt}}
\multiput(838.00,487.17)(8.509,-4.000){2}{\rule{0.600pt}{0.400pt}}
\put(849,482.17){\rule{2.500pt}{0.400pt}}
\multiput(849.00,483.17)(6.811,-2.000){2}{\rule{1.250pt}{0.400pt}}
\multiput(861.00,480.94)(1.505,-0.468){5}{\rule{1.200pt}{0.113pt}}
\multiput(861.00,481.17)(8.509,-4.000){2}{\rule{0.600pt}{0.400pt}}
\multiput(872.00,476.93)(0.611,-0.489){15}{\rule{0.589pt}{0.118pt}}
\multiput(872.00,477.17)(9.778,-9.000){2}{\rule{0.294pt}{0.400pt}}
\put(491.0,166.0){\rule[-0.200pt]{2.650pt}{0.400pt}}
\put(211.0,131.0){\rule[-0.200pt]{0.400pt}{100.696pt}}
\put(211.0,131.0){\rule[-0.200pt]{161.885pt}{0.400pt}}
\put(883.0,131.0){\rule[-0.200pt]{0.400pt}{100.696pt}}
\put(211.0,549.0){\rule[-0.200pt]{161.885pt}{0.400pt}}
\end{picture}
}
\subfloat[][Correction]{
\input{3d/m128m4_corr_p.tex}
}
\caption{Volume and error and related correction on a mesh with $128^3$ linear NURBS elements.}
\label{fig:3Dvol}
\end{center}
\end{figure}
It can be seen in figure \ref{fig:3Dvol} the volume is still quite well conserved, although not 
to full machine precision. To achieve this (near) conservation the required perturbation per 
timestep is still quite modest, similar as in the 2D case. 
 
\section{Conclusions}
\label{sec:conc}
In this paper four approaches to obtain a distance field in terms of element lengths have been presented. 
These approaches attempt to address two 
issues at the same time. To get truly smooth regularized Heaviside functions on arbitrary meshes, 
and to get a distance field without solving the difficult Eikonal equations directly.

In order to get smooth Heaviside functions on arbitrary meshes the  standard Eikonal equation has 
been recast in terms of a mesh-metric. Hereby the local meshsize has been taken into account
 when computing a distance field in terms of element lengths. 
Solving the Eikonal equation is circumvented by introducing a multiplicative scaling of the level-set. 
In this paper four different approaches for performing this scaling are investigated. For linear triangular meshes, linear quadrilateral meshes and $C_1$-continuous quadratic meshes the following conclusions can be drawn:
\begin{itemize}
\item {\it Direct redistancing}. The original level-set function is corrected using local derivative information. 
This results in a simple relation for the scaled distance field. 
 However, this approach fails to deliver a smooth Heaviside, even on $C_1$ meshes. Therefore, despite its simplicity, this approach has to be discarded. 
 
\item {\it Projected redistancing}. The discontinuous scaled distance field from the {\it direct redistancing} approach is projected on a continuous discretization. This corrects the smoothness issues but introduces other problems. 
On  distorted mesh the use of smoothing
results in spuriously diverging contour lines at the boundary. 
Additionally, in the vortex problem  insufficient smoothing results in a irregular solution.  
However increasing the smoothing leads to excessive problems to conserve the interface location. Therefore, this approach has to be 
discarded as well.

\item {\it Projected scaling}. Here only the scaling is projected instead of the scaled distance field itself.
This solves the problem of excessive interface movement at higher smoothing numbers. 
Unfortunately, low smoothing can still result in irregular solutions. We 
conclude that this is in principle a viable approach, but not the most desirable.

\item {\it Projected inverse scaling}. Here the reciprocal value is projected instead of the scaling itself.
In this approach the presence of smoothing plays a less prominent roll. 
In the absence of smoothing the approach still gives high quality solutions. 
Similar to  {\it projected scaling}  the interface  is precluded from moving. 
We conclude that this is approach has all the desired properties and is therefore the preferred choice.
\end{itemize}

A convergence study has been performed for {\it Projected inverse scaling}.
A $L_1$ 
convergence  study of the regularized Heaviside function has been done. 
This is effectively the mismatched area. Additionally, we have considered at the $L_{\infty}$ norm of the level-set.
It has been shown that the method is first order  accurate  when using linear triangles and quadrilaterals.
On quadratic quadrilaterals the method has a slightly higher order of accuracy, approaching an order of $1.5$.
This slightly reduced order of accuracy might be due to the global volume conservation. 

In the future the presented level-set approach will be employed to solve  two-fluid flow  potentially 
including (flexible) floating  objects.
By avoiding a solution of the redistancing problem, 
the new approaches should facilitate  strong coupling of the entire problem.
This lays the ground for strict energy control at the two-fluid interface.

\section*{Acknowledgements}

The author gratefully acknowledges  support from Delft University of Technology and Durham University.

\section*{References}
\bibliographystyle{unsrt}
\bibliography{references}


\end{document}